\documentclass{amsart}
\usepackage{amsmath,amsthm,amssymb}
\usepackage[mathcal]{euscript}

\usepackage[all]{xy}
\newcommand{\smallxymatrix}[1]{\xymatrix@1@=1pc
                                       {#1}
                                   } %for smaller diagrams
\newcommand{\tinyxymatrix}[1]{
\def\objectstyle{\scriptstyle}
\def\labelstyle{\scriptstyle}
\vcenter{\xymatrix@-1.2pc{#1} }} %for tiny diagrams

\newcommand{\inpair}[2]{\ar@<.5ex>[r]^-{#1}\ar@<-.5ex>[r]_-{#2} }
\newcommand{\smallpair}[4]{\smallxymatrix{
                      #1\ar@<.5ex>[r]^-{#2}\ar@<-.5ex>[r]_-{#3} & #4
                                         }} %for small parallel pair of arrows
\newcommand{\pair}[4]{\xymatrix@1{
                      #1\ar@<.5ex>[r]^-{#2}\ar@<-.5ex>[r]_-{#3} & #4
                                      }} %for parallel pair of arrows
\newcommand{\xylabel}[2]{\ar@{}[#1]|{#2}} %for tags in diagrams
\newcommand{\xylabelc}[2]{\ar@{}[#1]|-{#2}} %for tags in diagrams, centered
\newcommand{\xyinc}{\ar@{^{(}->}} %for incs in diagrams
\newcommand{\ten}{\mbox{\hspace*{-.5pt}\raisebox{1pt}{${\scriptstyle \otimes}$}
\hspace*{-4pt}}}
\newcommand{\cirprod}{\mbox{\hspace*{-.3pt}\raisebox{1pt}{${\scriptscriptstyle
\bigcirc}$}
\hspace*{-3.7pt}}}
\newcommand{\dirsum}{\mbox{\hspace*{-.5pt}\raisebox{1pt}{${\scriptstyle \oplus}$}
\hspace*{-3.5pt}}}

\newcommand{\tie}{\!\bowtie\!}
\newcommand{\ch}{{\textsf{char}}}      % FOR char(field)
\newcommand{\liebrac}[1]{[#1]}   % FOR LIE BRACKET
\newcommand{\eliebrac}{\liebrac{\ ,\ } }  % FOR [ , ]
\newcommand{\brac}[1]{\langle#1\rangle}   % FOR BRACE BRACKET
\newcommand{\mat}[1]{\begin{bmatrix}#1\end{bmatrix}} % for a matrix [ ]

\newcommand{\odash}[2]{\mbox{$#1$\nobreakdash-\hspace{0pt}#2}} %for A-module or n-dimensional
\newcommand{\infbi}{\odash{\epsilon}{bialgebra}}
\newcommand{\infbis}{\odash{\epsilon}{bialgebras}}
\newcommand{\infmod}{\odash{\epsilon}{Hopf module}}
\newcommand{\infmods}{\odash{\epsilon}{Hopf modules}}
\newcommand{\infbimod}{\odash{\epsilon}{Hopf bimodule}}
\newcommand{\infbimods}{\odash{\epsilon}{Hopf bimodules}}

\newcommand{\infbisubs}{\odash{\epsilon}{subbialgebras}}

\newcommand{\dendri}{\mbox{Dendriform algebras}}
\newcommand{\prelie}{\mbox{Pre-Lie algebras}}
\newcommand{\lie}{\mbox{Lie algebras}}
\newcommand{\assoc}{\mbox{Associative algebras}}
 \newcommand{\bracealg}{\mbox{Brace algebras}}

\newcommand{\teninf}{\ten_\epsilon}

\newcommand{\id}{\mathit{id}}

\newcommand{\gl}{\textsf{gl}}     %for lie algebra gl(V)

\newcommand{\Spa}{\textsf{Vec}}
\newcommand{\Alg}{\textsf{Alg}}
\newcommand{\AAlg}{\textsf{AAlg}}

\newcommand{\End}{\textsf{End}}

\newcommand{\Der}{\textsf{Der}}

\newcommand{\equal}[1]{
{\stackrel{{\textstyle #1}}{\ {\textstyle =}\ } }}

\newcommand{\map}[1]{\xrightarrow{#1}}

\newcommand{\inc}{\hookrightarrow}
\newcommand{\Z}{\mathbb{Z}}

\newcommand{\bfx}{\mathbf{x}}
\newcommand{\bfy}{\mathbf{y}}

\theoremstyle{plain}
\newtheorem{theo}{Theorem}[section]
\newtheorem{prop}[theo]{Proposition}
\newtheorem{lemm}[theo]{Lemma}
\newtheorem{coro}[theo]{Corollary}

\theoremstyle{definition}
\newtheorem{defi}[theo]{Definition}
\newtheorem{exas}[theo]{Examples} %use \bex
\newtheorem{exa}[theo]{Example}   %use \begin{exa}

\theoremstyle{remark}
\newtheorem{rems}[theo]{Remarks}  %use \br
\newtheorem{rem}[theo]{Remark}   %use \begin{rem}
\numberwithin{equation}{section}

\newcommand{\bd}{\begin{defi}}
\newcommand{\ed}{\end{defi}}
\newcommand{\bt}{\begin{theo}}
\newcommand{\et}{\end{theo}}
\newcommand{\bp}{\begin{prop}}
\newcommand{\ep}{\end{prop}}
\newcommand{\bl}{\begin{lemm}}
\newcommand{\el}{\end{lemm}}
\newcommand{\bc}{\begin{coro}}
\newcommand{\ec}{\end{coro}}
\newcommand{\bpf}{\begin{proof}}
\newcommand{\epf}{\end{proof}}
\newcommand{\br}{\begin{rems}\begin{flushleft}\end{flushleft}\nopagebreak} %remarks
\newcommand{\er}{\end{rems}}
\newcommand{\bex}{\begin{exas}\begin{flushleft}\end{flushleft}\nopagebreak} %examples
\newcommand{\eex}{\end{exas}}

\newcommand{\be}{\begin{enumerate}}
\newcommand{\ee}{\end{enumerate}}
\newcommand{\bi}{\begin{itemize}}
\newcommand{\ei}{\end{itemize}}

\begin{document}

\title[Infinitesimal, pre-Lie and dendriform]{Infinitesimal bialgebras, pre-Lie and
dendriform algebras}
\author[M. Aguiar]{Marcelo Aguiar}
\address{Department of Mathematics\\
         Texas A\&M University\\
         College Station\\
         Texas \ 77843\\
         USA}
\email{maguiar@math.tamu.edu}
\urladdr{http://www.math.tamu.edu/\~{}maguiar}
%\thanks{}

\date{October 31, 2002}
%16W30 Hopf algebras, 17A30 Algebras satisfying other identities
%17A42 Other $n$-ary compositions $(n \ge 3)$
%18D50-operads,  17D25-Lie admissible algebras.
\begin{abstract} We introduce the categories of infinitesimal Hopf modules and
bimodules over an infinitesimal bialgebra. We show that they
correspond to modules and bimodules over the infinitesimal version of the double.
We show that there is a natural, but non-obvious way to construct a pre-Lie
algebra from an arbitrary
infinitesimal bialgebra and a dendriform algebra from a quasitriangular infinitesimal
bialgebra. As consequences, we obtain a pre-Lie structure on the space of paths on
an arbitrary quiver, and a striking dendriform structure on the space of endomorphisms of
an arbitrary infinitesimal bialgebra, which combines the convolution and composition
products. We extend the previous constructions to the categories of Hopf, pre-Lie and
dendriform bimodules. We construct a brace algebra structure from an arbitrary
infinitesimal bialgebra;
this refines the pre-Lie algebra construction. In two appendices, we show that infinitesimal bialgebras are
comonoid objects in a certain monoidal category and discuss a related construction
for counital infinitesimal bialgebras.
\end{abstract}
\maketitle

\section{Introduction} \label{S:int}

The main results of this paper establish connections between
infinitesimal bialgebras, pre-Lie algebras and dendriform algebras, which were a priori
unexpected.

An infinitesimal bialgebra (abbreviated $\infbi$) is
a triple $(A, \mu,\Delta)$ where
$(A, \mu)$  is an  algebra,
$(A,\Delta)$  is a coalgebra, and $\Delta$ is a derivation (see Section~\ref{S:infmod}).
We write $\Delta(a)=a_1\ten a_2$, omitting the sum symbol.

Infinitesimal bialgebras were introduced by Joni and Rota~\cite[Section XII]{JR}. The
basic theory of these objects was developed  in~\cite{cor,analog}, where analogies with the
theories of ordinary Hopf algebras and Lie bialgebras were found; among which we remark the
existence of a ``double'' construction analogous to that of Drinfeld for
ordinary Hopf algebras or Lie bialgebras.
On the other hand, infinitesimal bialgebras have found important applications in
combinatorics~\cite{cdindex,ER}.

A pre-Lie algebra is a vector space $P$ equipped with an operation $x\circ y$   
satisfying a certain axiom~\eqref{E:prelie}, which guarantees that $x\circ y-y\circ x$
defines a Lie algebra structure on $P$.
These objects were introduced  by Gerstenhaber~\cite{Ger}, whose terminology
we follow, and independently by Vinberg~\cite{Vin}. See~\cite{Dzh,CL} for more references, examples,
and some of the general theory of pre-Lie algebras.

We show that any $\infbi$ can be turned into a pre-Lie algebra by defining 
\[a\circ b= b_1ab_2\,.\]
This is Theorem~\ref{T:inflie}. As an application, we construct a canonical pre-Lie
structure on the space of paths on an arbitrary quiver. We also note that the Witt
Lie algebra arises in this way from the
$\infbi$ of divided differences (Examples~\ref{Ex:prelie}). Other properties of this
construction are provided in Section~\ref{S:prelie}.

A dendriform algebra is a space $D$ equipped with two operations $x\succ y$ and $x\prec y$
 satisfying certain axioms~\eqref{E:dendri}, which guarantee that
$x\succ y+x\prec y$ defines an associative algebra structure on $D$. Dendriform algebras 
were introduced by Loday~\cite[Chapter 5]{L}. See~\cite{C,R,LR1,LR2} for additional recent
work on this subject.

There is a special class of $\infbis$ for which the derivation $\Delta$ is principal,
called quasitriangular $\infbis$. These are defined from solutions $r=\sum u_i\ten v_i$ of
the associative Yang-Baxter equation, introduced in~\cite{cor} and reviewed in
Section~\ref{S:infmod} of this paper. In Theorem~\ref{T:quasi-dendri}, we show that any
quasitriangular $\infbi$ can be made into a dendriform algebra by defining
\[x\succ y=\sum_i u_ixv_iy \text{ \ and \ }x\prec y=\sum_i xu_iyv_i\,.\]
This is derived from a more general construction of dendriform algebras from 
associative algebras equipped with a Baxter operator, given
in Proposition~\ref{P:baxter-dendri}. (Baxter operators should not be confused with Yang-Baxter
operators, see Remark~\ref{R:baxter}.)

As a main application of this construction, we work out the dendriform algebra
structure associated to the Drinfeld double of an $\infbi$ $A$.
This construction, introduced in~\cite{cor} and reviewed here in Section~\ref{S:infmod},
produces a quasitriangular $\infbi$ structure on the space $(A\ten A^*)\dirsum A\dirsum
A^*$. We provide explicit formulas for the resulting dendriform
structure in Theorem~\ref{T:double-dendri}. This is one of the main results of
this paper. It turns out that the subspace
$A\ten A^*$ is closed under the dendriform operations. The resulting dendriform algebra
structure on the space $\End(A)$ of linear endomorphisms of $A$ is
(Corollary~\ref{C:end-dendri})
\[T\succ S=(\id\ast T\ast id)S +(\id\ast T)(S\ast \id) \text{ \ and \ }
T\prec S=T(\id\ast S\ast id)+(T\ast \id)(\id\ast S)\,.\]
In this formula, $T$ and $S$ are arbitrary endomorphisms of $A$, $T\ast S=\mu(T\ten
S)\Delta$ is the convolution, and the concatenation of endomorphisms denotes composition.
 When $A$ is a quasitriangular
$\infbi$, our results give dendriform structures on $A$ and $\End(A)$. In Proposition~\ref{P:pi-dendri},
we show that they are related by a canonical morphism of dendriform algebras $\End(A)\to A$.

Other properties of the construction of dendriform algebras  are given in 
Section~\ref{S:dendri}. In particular, it is shown that the
constructions of pre-Lie algebras from $\infbis$ and of dendriform algebras from quasitriangular
$\infbis$ are compatible, in the sense that the diagram
\[\xymatrix@1{ {\mbox{Quasitriangular
$\infbis$}}\ar[d]\ar[r] & {\infbis}\ar[d]\\ {\dendri}\ar[r] & {\prelie}  } \]
commutes.

This paper also introduces the appropriate notion of modules over infinitesimal bialgebras. 
These are called infinitesimal Hopf modules, abbreviated $\infbimods$. They are defined in
Section~\ref{S:infmod}. In the same section, it is shown that $\infbimods$  are
precisely modules over the double, when the
$\infbi$ is finite dimensional (Theorem~\ref{T:moddouble}), and that any module can be turned
into an $\infbimod$, when the $\infbi$ is quasitriangular
(Proposition~\ref{P:modquasi}).

The constructions of dendriform and pre-Lie algebras are extended to the corresponding
categories of bimodules in Section~\ref{S:bimod}. A commutative diagram of the form 
\[\xymatrix{ {\mbox{Associative bimodules}}\ar[d]\ar[r] & {\infbimods}\ar[d]\\
{\mbox{Dendriform bimodules}}\ar[r] & {\mbox{Pre-Lie bimodules}}  } \] is obtained. 

A brace algebra is a space $B$ equipped with a family of higher degree operations satisfying
certain axioms\eqref{E:brace}. Brace algebras originated in work of Kadeishvili~\cite{Kad},
Getzler~\cite{Get} and Gerstenhaber and Voronov~\cite{GV95a,GV95b}. In this paper we deal
with the ungraded, unsigned version of these objects, as in the recent works of
Chapoton~\cite{C} and Ronco~\cite{R}. Brace algebras sit between dendriform and
pre-Lie; as explained in~\cite{C,R}, the functor from dendriform to  pre-Lie algebras factors through the category
of brace algebras. Following a suggestion of Ronco, we show in Section~\ref{S:brace} that the
construction of pre-Lie algebras from $\infbis$ can be refined accordingly. We associate a
brace algebra to any
$\infbi$ (Theorem~\ref{T:infbrace}) and obtain a commutative diagram
\[\xymatrix@1{ {\mbox{Quasitriangular
$\infbis$}}\ar[d]\ar[r] & {\infbis}\ar[d]\ar[rd]\\ {\dendri}\ar[r] & {\bracealg}\ar[r] &
{\prelie}  }\]
The brace algebra associated to the $\infbi$ of divided differences is explicitly described
in Example~\ref{Ex:brace}. The higher braces are given by
\[\brac{\bfx^{p_1},\ldots,\bfx^{p_n};\bfx^r}=\binom{r}{n}\bfx^{r+p_1+\cdots+p_n-n}\,, \]
where $\binom{r}{n}$ is the binomial coefficient.

In Appendix~\ref{S:comonoid} we construct a certain
monoidal category of algebras for which the comonoid objects are precisely $\infbis$, and we
discuss how $\infbis$ differ from bimonoid objects in certain related braided monoidal
categories.

In Appendix~\ref{S:counital} we study certain special features of 
counital $\infbis$.  We construct another monoidal category of algebras and show that
comonoid objects in this category are precisely counital
$\infbis$ (Proposition~\ref{P:comonoid}). The relation to the constructions of
Appendix~\ref{S:comonoid} is explained.  We also describe counital
$\infmods$ in terms of this monoidal structure (Proposition~\ref{P:infmod}).

\subsection*{Notation and basic terminology} All spaces and algebras are over a fixed field
$k$, often omitted from the notation. Sum symbols are omitted from Sweedler's notation: we write
$\Delta(a)=a_1\ten a_2$ when
$\Delta$ is a coassociative comultiplication, and similarly for comodule structures. The
composition of maps $f:U\to V$ with $g:V\to W$ is denoted by $gf:U\to W$.

\section{Infinitesimal modules over infinitesimal bialgebras} \label{S:infmod}

An infinitesimal bialgebra (abbreviated $\infbi$) is 
a triple $(A, \mu,\Delta)$ where
$(A, \mu)$  is an  algebra,  $(A,\Delta)$  is a coalgebra,
and for each $a,b\in A$,
\begin{equation}
\Delta(ab)= ab_1\ten b_2+ a_1\ten a_2b\, .\label{E:infbi}
\end{equation} 
We do not require the algebra to be unital or the coalgebra to be counital.

A derivation of an algebra $A$ with values in a \odash{A}{bimodule} $M$  is a
linear map $D:A\to M$ such that
\[D(ab)=a\cdot D(b)+D(a)\cdot b\ \ \forall\ a,b\in A\ .\]
We view $A\ten A$ as an \odash{A}{bimodule} via
\[a\cdot(b\ten c)=ab\ten c \text{ and }(b\ten c)\cdot a=b\ten ca\,.\]

A coderivation from a \odash{C}{bicomodule} $M$ to a coalgebra $C$ is a map
$D:M\to C$ such that 
\[\Delta D=(\id_C\ten D) t+(D\ten\id_C) s\ , \]
where $t:M\to C\ten M$ and $s:M\to M\ten C$ are the bicomodule structure maps [Doi].
We view $C\ten C$ as a \odash{C}{bicomodule} via
\[t=\Delta\ten\id_C \text{ and }s=\id_C\ten\Delta\,.\]

The compatibility condition \eqref{E:infbi} may be written as 
\[\Delta \mu=( \mu\ten id_A)( id_A\ten\Delta)+( id_A\ten \mu)(\Delta\ten id_A
)\]
This says that $\Delta:A\to A\ten A$ is a derivation of 
the algebra $(A, \mu)$ with values in the
\odash{A}{bimodule} $A\ten A$, or equivalently, that
$ \mu:A\ten A\to A$ is a coderivation  from the
\odash{A}{bicomodule} $A\ten A$ with values in the coalgebra 
$(A,\Delta)$.

\bd\label{D:infmod} Let $(A,\mu,\Delta)$ be an $\infbi$.
A left infinitesimal Hopf module (abbreviated $\odash{\epsilon}{Hopf}$ module) over $A$ is 
a space $M$ endowed with a left \odash{A}{module} structure
$\lambda:A\ten M\to M$  and a left \odash{A}{comodule} structure
$\Lambda:M\to A\ten M$,  such that
\[\Lambda\lambda=( \mu\ten id_M)( id_A\ten\Lambda)+( id_A\ten\lambda)(\Delta\ten id_M)\,.\]
\ed
We will often write 
\[\lambda(a\ten m)=am \text{ and }\Lambda(m)=m_{-1}\ten m_0\,\]
The compatibility condition above may be written as  
$\Lambda(am)=a\Lambda(m)+\Delta(a)m$, or more explicitly,
\begin{equation}
(am)_{-1}\ten(am)_{0}= am_{-1}\ten m_0+ a_1\ten a_2m\,, \text{ \ for each $a\in
A$ and $m\in M$}.\label{E:infmod}
\end{equation}

The notion of $\infmods$ bears a certain analogy to the notion of Hopf modules over
ordinary Hopf algebras. The basic examples of Hopf modules from~\cite[1.9.2-3]{mon} admit
the following versions in the context of $\infbis$.
\begin{exas}\label{Ex:nontrivial} Let $(A,\mu,\Delta)$ be an $\infbi$.
\be 
\item $A$ itself is an $\infmod$ via $\mu$ and $\Delta$, precisely by definition
of $\infbi$.
\item More generally, for any space $V$, $A\ten V$ is an $\infmod$ via
\[\mu\ten\id:A\ten A\ten V\to A\ten V \text{ and }
\Delta\ten\id:A\ten V\to A\ten A\ten V\,.\]
\item A more interesting example follows. Assume that the coalgebra $(A,\Delta)$ admits a
counit
$\eta:A\to k$. Let
$N$ be a left \odash{A}{module}.   Then there is an $\infmod$ structure on the space
$A\ten N$ defined by
\[a\cdot(a'\ten n)=aa'\ten n+\eta(a')\,a_1\ten a_2n \text{ \ and \ }
\Lambda(a\ten n)=a_1\ten a_2\ten n\,.\]
This can be checked by direct calculations. A more conceptual proof will be given
later (Corollary~\ref{C:nontrivial}). Note that if $N$ is a trivial \odash{A}{module}
($an\equiv 0$) then this structure reduces to that of example 2.
\ee
\end{exas}

When $H$ is a finite dimensional (ordinary) Hopf algebra, left Hopf modules over $H$ are
precisely left modules over the {\em Heisenberg double} of $H$~\cite[Examples
4.1.10 and 8.5.2]{mon}. 

There is an analogous result for infinitesimal bialgebras which, as it turns out,
involves the  Drinfeld double of $\infbis$. 
 
 We first recall the construction of the Drinfeld double $D(A)$ of a finite
dimensional $\infbi$ $(A,\mu,\Delta)$ from~\cite[Section 7]{cor}. 
Consider the following version of the dual of $A$
\[A':=(A^*,\Delta^{*^{op}},-\mu^{*^{cop}})\ .\]
Explicitly, the structure on $A'$ is:
\begin{gather}
\label{E:proddual}
(f\cdot g)(a)=g(a_1)f(a_2)\ \forall\ a\in A,\ f,g\in A' \text{
and}\\
\label{E:coproddual} \Delta(f)=f_1\ten f_2 \iff f(ab)=-f_2(a)f_1(b)\ \forall\ f\in
A',\ a,b\in A\ .
\end{gather}
Below we always refer to this structure when dealing with multiplications or
comultiplications of elements of $A'$. Consider also the actions of $A'$ on $A$ and
$A$ on $A'$ defined by
\begin{gather}
f\to a=f(a_1)a_2 \text{ and } f\gets a =-f_2(a)f_1 \label{E:actionone}
\intertext{or equivalently}
g(f\to a)=(gf)(a) \text{ and }(f\gets a)(b)=f(ab)\ .\label{E:actiontwo}
\end{gather}

\bp \label{P:double}
Let $A$ be a finite dimensional $\infbi$, consider the vector space
\[D(A):=(A\ten A')\dirsum A\dirsum A'\]
and denote the element $a\ten f\in A\ten A'\subseteq D(A)$ by $a\tie f$.
Then $D(A)$ admits a unique $\infbi$ structure such that:
\bi
\item[(a)] $A$ and $A'$ are subalgebras, $a\cdot f=a\tie f$, $f\cdot a=f\to
a+f\gets a$, and
\item[(b)] $A$ and $A'$ are subcoalgebras. 
\ei
\ep
\bpf See~\cite[Theorem 7.3]{analog}. \epf

We will make use of the following universal property of the double.
\bp \label{P:univdouble}
Let $A$ be a finite dimensional $\infbi$, $B$ an algebra and $\rho:A\to B$ and $\rho':A'\to
B$ morphisms of algebras such that $\forall\ a\in A,\ f\in A'$,
\begin{equation} \label{E:univdouble}
\rho'(f)\rho(a)=\rho(f\to a)+\rho'(f\gets a)\ .
\end{equation}
Then there exists a unique morphism of algebras $\hat{\rho}:D(A)\to B$ such that
$\hat{\rho}_{|_A}=\rho$ and $\hat{\rho}_{|_A'}=\rho'$.
\ep
\bpf This follows from Propositions 6.5 and 7.1 in~\cite{cor}.
\epf

We can now show that $\infmods$ are precisely modules over the double.

\bt\label{T:moddouble}
Let $A$ be a finite dimensional $\infbi$ and $M$ a space.
 If $M$ is a left $\infmod$ over $A$ via $\lambda(a\ten m)=am$ and
$\Lambda(m)=m_{-1}\ten m_0$, then $M$ is a left module over $D(A)$ via
\[a\cdot m= am,\ \ f\cdot m=f(m_{-1})m_0 \text{ \ and \ }(a\tie f)\cdot
m=f(m_{-1})am_0\,.\]
 Conversely, if $M$ is a left module over $D(A)$, then $M$ is a left $\infmod$ over $A$
and the structures are related as above.
\et
\bpf Suppose first that $M$ is a left $\infmod$ over $A$. 

Since $(M,\Lambda)$ is a left \odash{A}{comodule}, it is also a left \odash{A'}{module}
via $f\cdot m:=f(m_{-1})m_0$. Let $\rho:A\to\End(M)$ and $\rho':A'\to\End(M)$ be the
morphisms of algebras corresponding to the left module structures:
\[\rho(a)(m)=am,\ \ \rho'(f)(m)=f(m_{-1})m_0\,.\]
We will apply Proposition~\ref{P:univdouble} to deduce the existence of a morphism
of algebras $\hat{\rho}:D(A)\to\End(M)$ extending $\rho$ and $\rho'$. We need to check
\eqref{E:univdouble}. We have
\begin{align*}
\rho'(f)\rho(a)(m)&=f\bigl((am)_{-1}\bigr)(am)_0\equal{\eqref{E:infmod}}
f(a_1)a_2m+f(am_{-1})m_0 \\
&\equal{(\ref{E:actionone},\ref{E:actiontwo})} (f\to a)m+(f\gets
a)(m_{-1})m_0\\ 
&=\rho(f\to a)(m)+\rho'(f\gets a)(m)
\end{align*}
as needed. Thus, $\hat{\rho}$ exists and $M$ becomes a left \odash{D(A)}{module}
via $\alpha\cdot m=\hat{\rho}(\alpha)(m)$. Since $\hat{\rho}$ extends $\rho$ and
$\rho'$, we have 
\begin{gather*}
a\cdot m=\rho(a)(m)=am,\ \ f\cdot m=\rho'(f)(m)=f(m_{-1})m_0\\
\intertext{and, from the description of the multiplication in $D(A)$ in
Proposition~\ref{P:double},} 
(a\tie f)\cdot m=\hat{\rho}(af)(m)=\rho(a)\rho'(f)(m)=f(m_{-1})am_0\,.
\end{gather*}
This completes the proof of the first assertion.

Conversely, if $M$ is a left \odash{D(A)}{module}, then restricting via the morphisms
of algebras $A\inc D(A)$ and $A'\inc D(A)$, $M$ becomes a left \odash{A}{module}
and left \odash{A'}{module}. As above, the latter structure is equivalent to a left
\odash{A}{comodule} structure on $M$. From the associativity axiom 
\[f\cdot(a\cdot m)=(fa)\cdot m=(f\to a)\cdot m+(f\gets a)\cdot m\]
we deduce
\[f\bigl((am)_{-1}\bigr)(am)_0=f(a_1)a_2m+f(am_{-1})m_0\,.\]
Since this holds for every $f\in A'$, we obtain  the $\infmod$
Axiom~\eqref{E:infmod}. Also,
\[(a\tie f)\cdot m=(af)\cdot m=a\cdot(f\cdot m)=f(m_{-1})am_0\,,\]
so the structures of left module over $D(A)$ and left $\infmod$ over $A$ are related
as stated.
\epf

We close the section by showing that when $A$ is a quasitriangular $\infbi$,
any  \odash{A}{module} carries a natural structure of $\infmod$ over $A$.

We first recall the definition of quasitriangular $\infbis$.
Let $A$ be an associative
algebra. An element $r=\sum_i u_i\ten v_i\in A\ten A$ is a solution of the associative
Yang-Baxter equation~\cite[Section 5]{cor} if
\begin{gather}
r_{13}r_{12}-r_{12}r_{23}+r_{23}r_{13}=0 \label{E:ayb}\\
\intertext{or, more explicitly,}
\sum_{i,j} u_iu_j\ten v_j\ten v_i-\sum_{i,j} u_i\ten v_iu_j\ten
v_j+\sum_{i,j} u_j\ten u_i\ten v_iv_j=0\,.\notag
\end{gather}
This condition implies that the {\em principal} derivation $\Delta:A\to A\ten A$ defined by
\begin{equation}\label{E:principal}
\Delta(a)=r\cdot a-a\cdot r=\sum_i u_i\ten v_ia-\sum_i au_i\ten v_i\,,
\end{equation}  
is coassociative~\cite[Proposition 5.1]{cor}. 
Thus, endowed with this comultiplication, $A$
becomes an $\infbi$. We refer to the pair $(A,r)$ as a quasitriangular 
$\infbi$~\cite[Definition 5.3]{cor}.

\begin{rem}\label{R:signs} In our previous work~\cite{cor,analog}, we have used the
comultiplication
\[-\Delta(a)=a\cdot r-r\cdot a=\sum_i au_i\ten v_i-\sum_i u_i\ten v_ia \] 
instead of $\Delta$. Both $\Delta$ and $-\Delta$ endow $A$ with a structure of $\infbi$,
and there is no essential difference in working with one or the other. The choice we
adopt in ~\eqref{E:principal}, however, is more
convenient for the purposes of this work, particularly in relating
quasitriangular $\infbis$ and their bimodules to dendriform algebras and their bimodules
(Sections~\ref{S:dendri} and~\ref{S:bimod}).

It is then necessary to make the corresponding sign adjustments to the results on
quasitriangular $\infbis$ from~\cite{cor,analog} before applying them in the present
context. For instance, Proposition 5.5 in~\cite{cor} translates as
\begin{equation}\label{E:deltaonr}
\sum_{i,j} u_i\ten u_j\ten v_jv_i=
r_{23}r_{13}=(\Delta\ten\id)(r)=\Delta(u_i)\ten v_i\,.
\end{equation}
\end{rem}

\bp\label{P:modquasi} Let $(A,r)$ be a quasitriangular $\infbi$ and $M$ a left
\odash{A}{module}. Then $M$ becomes a left $\infmod$ over $A$ via $\Lambda:M\to A\ten M$,
\[\Lambda(m)=\sum_i u_i\ten v_im\,.\]
\ep
\bpf We first check that $\Lambda$ is coassociative, i.e., 
$(\id\ten\Lambda)\Lambda=(\Delta\ten\id)\Lambda$. We have
\begin{gather*}
(\id\ten\Lambda)\Lambda(m)=\sum_i u_i\ten\Lambda(v_im)=\sum_{i,j} u_i\ten u_j\ten
v_jv_im\\
\intertext{and}
(\Delta\ten\id)\Lambda(m)=\sum\Delta(u_i)\ten v_im\,.
\end{gather*}
According to~\eqref{E:deltaonr}, these two expressions agree.

It only remains to check Axiom~\eqref{E:infmod}. Since
$\Delta(a)=\sum_i u_i\ten v_ia-\sum_i au_i\ten v_i$,
we have
\[\Delta(a)m+a\Lambda(m)=\sum_i u_i\ten v_iam-\sum_i au_i\ten v_im+\sum_i au_i\ten
v_im=\sum_i u_i\ten v_iam=\Lambda(am)\,,\]
as needed.
\epf

\begin{rem} If $A$ is a finite dimensional quasitriangular $\infbi$, then there is a
canonical morphism of $\infbis$ $\pi:D(A)\to A$, which is the identity on
$A$~\cite[Proposition 7.5]{cor}. Therefore, any left \odash{A}{module} $M$ can be
first made into a left \odash{D(A)}{module} by restriction via $\pi$, and then, by
Theorem~\ref{T:moddouble}, into a left $\infmod$ over $A$. 
It is easily seen that this structure coincides with the one of
Proposition~\ref{P:modquasi}. Note that the construction of the latter proposition
is more general, since it does not require finite dimensionality of $A$.
\end{rem}

\section{Pre-Lie algebras} \label{S:prelie}

\bd \label{D:prelie} A (left) pre-Lie algebra is a vector space $P$ together with a
 map $\circ:P\ten P\to P$ such that
\begin{equation} \label{E:prelie}
x\circ(y\circ z)-(x\circ y)\circ z=y\circ(x\circ z)-(y\circ x)\circ z\,.
\end{equation}
\ed
 There is a similar notion of right pre-Lie algebras. In this paper, we will only deal with
left pre-Lie algebras and we will refer to them simply as pre-Lie algebras.

Defining a new operation $P\ten P\to P$ by $\liebrac{x,y}=x\circ y-y\circ x$
one obtains a  Lie algebra structure on $P$~\cite[Theorem 1]{Ger}.

Next we show that every $\infbi$ $A$ gives rise to a structure of pre-Lie algebra,
and hence also of Lie algebra, on the underlying space of $A$.

\bt \label{T:inflie} Let $(A,\mu,\Delta)$ be an $\infbi$. Define a new operation on $A$ by
\begin{equation}\label{E:inflie}
a\circ b= b_1ab_2\,.
\end{equation}
Then $(A,\circ)$ is a pre-Lie algebra.
\et
\bpf By repeated use of~\eqref{E:infbi} we find
\[\Delta(abc)=ab\cdot\Delta(c)+\Delta(ab)\cdot
c= abc_1\ten c_2+ab_1\ten b_2c+a_1\ten a_2bc\,.\]
Together with coassociativity this gives
\[\Delta(c_1bc_2)= c_1bc_2\ten c_3+c_1b_1\ten b_2c_2+c_1\ten c_2bc_3\,.\]
Combining this with~\eqref{E:inflie} we obtain
\[a\circ(b\circ c)=a\circ(c_1bc_2)=c_1bc_2a c_3+c_1b_1a b_2c_2+c_1a c_2bc_3\,.\]
On the other hand,
\[(a\circ b)\circ c=(b_1ab_2)\circ c=c_1b_1ab_2c_2\,.\]
Therefore,
\[a\circ(b\circ c)-(a\circ b)\circ c=c_1bc_2a c_3+c_1a c_2bc_3\,.\]
Since this expression is invariant under $a\leftrightarrow b$, Axiom \eqref{E:prelie} holds
and
$(A,\circ)$ is a pre-Lie algebra.
\epf

For a vector space $V$, let $\gl(V)$ denote the space of all linear maps
$V\to V$, viewed as a Lie algebra under the commutator bracket $[T,S]=TS-ST$.

If $P$ is a pre-Lie algebra and $x\in P$, let $L_x:P\to P$ be
$L_x(y)=x\circ y$. The map $L:P\to \gl(P)$,
$x\mapsto L_x$ is a morphism of Lie algebras. This statement is just a reformulation
of Axiom~\eqref{E:prelie}.

In the case when the pre-Lie algebra comes from an $\infbi$ $(A,\mu,\Delta)$, more can
be said about this canonical map. Let $\Der(A,\mu)$ denote the space of all
derivations $D:A\to A$ of the associative algebra $A$. Recall that this is a Lie
subalgebra of $\gl(A)$.

\bp\label{P:derivation} Let $(A,\mu,\Delta)$ be an $\infbi$ and consider the associated
pre-Lie and Lie algebra structures on $A$.
The canonical morphism of Lie
algebras $L:(A,\eliebrac)\to\gl(A)$ actually maps to the Lie subalgebra $\Der(A,\mu)$
of $\gl(A)$.
\ep
\bpf We must show that each $L_c\in\gl(A)$ is a derivation of
the associative algebra $A$. We have
$\Delta(ab)\equal{\eqref{E:infbi}}ab_1\ten b_2+a_1\ten a_2b$ and hence
\[L_c(ab)=c\circ(ab)\equal{\eqref{E:inflie}}ab_1cb_2+a_1ca_2b\equal{\eqref{E:inflie}}
a(c\circ b)+(c\circ a)b=aL_c(b)+L_c(a)b\ ,\]
as needed.

\epf

\bex\label{Ex:prelie}
\be
\item Consider the $\infbi$ of {\em divided differences}. This is the
algebra $k[\bfx,\bfx^{-1}]$ of Laurent polynomials, with
$\Delta(f(\bfx))=\frac{f(\bfx)-f(\bfy)}{\bfx-\bfy}$. This was the example that
motivated Joni and Rota to abstract the notion of  $\infbis$~\cite[Section XII]{JR}.
More explicitly,
\[\Delta(\bfx^n)=\sum_{i=0}^{n-1}\bfx^i\ten\bfx^{n-1-i} \text{ \ and \ }
\Delta(\frac{1}{\bfx^n})=-\sum_{i=1}^{n}\frac{1}{\bfx^i}\ten\frac{1}{\bfx^{n+1-i}}\,,
\text{ \ for $n\geq 0$.}\]
The corresponding pre-Lie algebra structure is
\[\bfx^m\circ\bfx^n=n\bfx^{m+n-1}\,,\text{ \ for any $n\in\Z$,} \]
and the Lie algebra structure on $k[\bfx,\bfx^{-1}]$ is
\[\liebrac{\bfx^m,\bfx^n}=(n-m)\bfx^{m+n-1} \text{ \ for  $n$, $m\in\Z$.} \]
This is the so called Witt Lie algebra. The canonical map
$k[\bfx,\bfx^{-1}]\to\Der(k[\bfx,\bfx^{-1}],\mu)$ of Proposition~\ref{P:derivation}
sends $\bfx^m$ to $\bfx^m\frac{d}{d\bfx}$, so it is an isomorphism of Lie algebras.

\item The algebra of matrices $M_2(k)$ is an $\infbi$  under
\[\Delta\mat{a & b\\c & d}=\mat{0 & a\\0 & c}\ten\mat{0 & 1\\0 & 0}-
\mat{0 & 1\\0 & 0}\ten\mat{c & d\\0 & 0}\]
\cite[Example 2.3.7]{cor}. One finds easily
that the corresponding Lie algebra splits as a direct sum of Lie algebras
\[\mathfrak{g}=\mathfrak{h}\dirsum \mathfrak{o}\]
where $\mathfrak{h}=k\{x,y,z\}$ is the \odash{3}{dimensional} {\em Heisenberg} algebra
\[\{x,y\}=z,\ \ \{x,z\}=\{y,z\}=0\ ,\]
and $\mathfrak{o}=k\{i\}$ is the \odash{1}{dimensional} Lie algebra. To
realize this isomorphism explicitly, one may take
\[x=\mat{0 & 0\\1 &0},\ y=\mat{1 & 0\\0 &0},\ z=\mat{0 & 1\\0 & 0} \text{ and }
i=\mat{1 & 0\\0 &1}\ .\]

\item  The {\em path algebra} of a quiver carries a canonical $\infbi$
structure~\cite[Example 2.3.2]{cor}.  Let
$Q$ be an arbitrary quiver (i.e., an oriented graph). Let $Q_n$ be the set of paths
$\alpha$ in $Q$ of length $n$:
\[\alpha: e_0\map{a_1}e_1\map{a_2}e_2\map{a_3}\ldots e_{n-1}\map{a_n}e_n\ .\]
In particular, $Q_0$ is the set of vertices and $Q_1$ is the set of arrows.
 Recall that the path algebra of $Q$ is the space
$kQ=\dirsum_{n=0}^\infty kQ_n$ where multiplication is concatenation of paths whenever
possible; otherwise is zero. The comultiplication  is defined on a path
$\alpha=a_1a_2\ldots a_n$ as above by
\[\Delta(\alpha)=e_0\ten a_2a_3\ldots a_n+a_1\ten a_3\ldots a_n+\ldots+
a_1\ldots a_{n-1}\ten e_n\ .\]
In particular, $\Delta(e)=0$ for every vertex $e$ and $\Delta(a)=e_0\ten
e_1$ for every arrow $e_0\map{a}e_1$.

In order to describe the corresponding pre-Lie algebra structure on $kQ$, consider pairs
$(\alpha,b)$ where $\alpha$ is a path from $e_0$ to $e_n$ (as above) and $b$ is an arrow
from $e_0$ to $e_n$. Let us call such a pair a {\em shortcut}. The  pre-Lie algebra
structure on $kQ$ is
\[\alpha\circ\beta=\sum_{b_i\in\beta} b_1\ldots b_{i-1}\alpha b_{i+1}\ldots b_m\,,\]
where the sum is over all arrows $b_i$ in the path $\beta=b_1\ldots b_m$ such that
$(\alpha,b_i)$ is a shortcut.
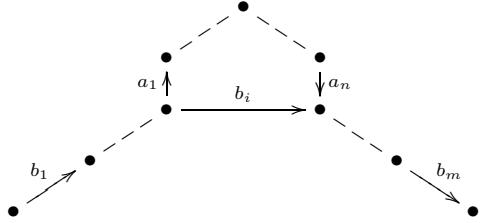
\[\xymatrix@C-6pt@R-15pt{
 & & & {\bullet}\ar@{--}[rd] & & &\\
& & {\bullet}\ar@{--}[ru]& & {\bullet}\ar[d]^{a_n} & &\\
& & {\bullet}\ar[u]^{a_1}\ar[rr]^{b_i}& &{\bullet}\ar@{--}[rd] & &\\
&{\bullet}\ar@{--}[ru]& & & &{\bullet}\ar[rd]^{b_m}&\\
{\bullet}\ar[ru]^{b_1} & & & & & & {\bullet}
}\]
\ee
\eex

\medskip

A {\em biderivation} of an $\infbi$ $(A,\mu,\Delta)$ is a map $B:A\to A$
that is both a derivation of $(A,\mu)$ and a coderivation of $(A,\Delta)$, i.e.,
\begin{equation}\label{E:bider}
B(ab)=aB(b)+B(a)b \text{ \ and \ }\Delta(B(a))=a_1\ten
B(a_2)+B(a_1)\ten a_2\,.
\end{equation}
 A derivation of a pre-Lie algebra is a map $D:P\to P$ such that
\[D(x\circ y)=x\circ D(y)+D(x)\circ y\,. \]
Such a map $D$ is always a derivation of the associated Lie algebra.

\bp\label{P:bider} Let $B$ be a biderivation of an $\infbi$ $A$. Then $B$ is a
derivation of the associated pre-Lie algebra (and hence also of the associated Lie
algebra).
\ep
\bpf We have
\[B(a)\circ b=b_1B(a)b_2 \text{ \ and \ } a\circ B(b)=b_1aB(b_2)+B(b_1)ab_2\,.\]
Hence,
\[B(a)\circ b+a\circ B(b)=b_1B(a)b_2+b_1aB(b_2)+B(b_1)ab_2=B(b_1ab_2)=B(a\circ b)\,.\]
\epf

The construction of a pre-Lie algebra from and $\infbi$ can be extended to the categories
of modules. This will be discussed in the appropriate generality in Section~\ref{S:bimod}.
A first result in this direction is discussed next.

Let $(P,\circ)$ be a pre-Lie algebra.
A  left \odash{P}{module} is a space
$M$ together with a map $P\ten M\to M$, $x\ten m\mapsto x\circ m$, such that
\begin{equation}\label{E:preliemod}
x\circ(y\circ m)-(x\circ y)\circ m=y\circ(x\circ m)-(y\circ x)\circ m\,.
\end{equation}

\bp\label{P:preliemod} Let $A$ be an $\infbi$ and $M$ a left $\infmod$ over $A$ via
\[\lambda(a\ten m)=am \text{ \ and \ }\Lambda(m)=m_{-1}\ten m_0\,.\]
Then $M$ is a left pre-Lie module over the pre-Lie algebra $(A,\circ)$ of
Theorem~\ref{T:inflie} via
\[a\circ m= m_{-1}am_0\,.\]
\ep
\bpf We first compute
\begin{align*}
\Lambda(a\circ m) &=\Lambda(m_{-1}am)\equal{\eqref{E:infmod}}
\Delta(m_{-1}a)m_0+m_{-1}a\Lambda(m_0)\\
&\equal{\eqref{E:infbi}}m_{-1}\ten m_{-1}am_0+m_{-1}a_1\ten a_2 m_0+m_{-2}am_{-1}\ten
m_0\,,
\end{align*}
where we have used the coassociativity axiom for the comodule structure $\Lambda$.
It follows that
\[b\circ(a\circ m)=m_{-2}bm_{-1}am_0+m_{-1}a_1b a_2 m_0+m_{-2}am_{-1}bm_0\,.\]
On the other hand,
\[(b\circ a)\circ m=m_{-1}(b\circ a)m_0\equal{\eqref{E:inflie}}m_{-1}a_1ba_2m_0\,.\]
Therefore,
\[b\circ(a\circ m)-(b\circ a)\circ m =m_{-2}bm_{-1}am_0+m_{-2}am_{-1}bm_0\,.\]
Since this expression is invariant under $a\leftrightarrow b$, Axiom~\eqref{E:preliemod}
holds.
\epf

\begin{rem}\label{R:dual}
Since the notion of $\infbis$ is self-dual~\cite[Section 2]{cor}, one should expect
dual constructions to those of Theorem~\ref{T:inflie} and  Propositions~\ref{P:bider}
and~\ref{P:preliemod}. This is indeed the case. Namely, if $A$ is an arbitrary $\infbi$,
then the map
$\gamma:A\to A\ten A$ defined by
\[\gamma(a)=a_2\ten a_1a_3\]
endows $A$ with a structure of {\em left pre-Lie coalgebra}. Also, if $B:A\to A$ is
a biderivation of $A$ then it is also a coderivation of $(A,\gamma)$. Moreover, if
$M$ is a left
$\infmod$ over $A$, then $M$ is left pre-Lie comodule over $(A,\gamma)$ via $\psi:M\to
A\ten M$ defined by
\[\psi(m)=m_{-1}\ten m_{-2}m_0\,.\]

If $A$ is an $\infbi$, then $A$ carries structures of pre-Lie
algebra and pre-Lie coalgebra, as just explained. Hence, it also
carries structures of Lie algebra and Lie coalgebra, by
\[\liebrac{a,b}=b_1ab_2-a_1ba_2 \text{ \ and \ }\delta(a)=a_2\ten a_1a_3-a_1a_3\ten
a_2\,.\]
 In general, these structures are not compatible, in the sense that they do {\em 
not} define a structure of Lie bialgebra on $A$. 
\end{rem}

\section{Dendriform algebras} \label{S:dendri}

\bd\label{D:dendri} A dendriform algebra is a vector space $D$
together with  maps
$\succ :D\ten D\to D$ and $\prec :D\times D\to D$ such that
\begin{align}
(x\prec y)\prec z &= x\prec (y\prec z)+x\prec (y\succ z) \notag\\
x\succ (y\prec z) &= (x\succ y)\prec z \label{E:dendri}\\
x\succ (y\succ z) &=(x\prec y)\succ z+(x\succ y)\succ z\,. \notag
\end{align}
\ed

Dendriform algebras were introduced by Loday~\cite[Chapter 5]{L}.
There is also a notion of {\em dendriform trialgebras}, which involves three
operations~\cite{LR3}. When it is necessary to distinguish between these two notions,
one uses the name {\em dendriform dialgebras} to refer to what in this paper (and
in~\cite{L}) are called dendriform algebras. Since only dendriform algebras (in the sense
of Definition~\ref{D:dendri}) will be considered in this paper, this usage will not
be adopted.

Let $(D,\succ,\prec)$ be a dendriform algebra. Defining  $x\cdot
y=x\succ y+x\prec y$ one obtains an associative algebra structure on $D$. In addition,
defining
\begin{equation}\label{E:dendprelie}
x\circ y=x\succ y-y\prec x
\end{equation}
 one obtains a (left) pre-Lie algebra structure on $D$.  Moreover, the Lie algebras
canonically associated to $(D,\cdot)$ and $(D,\circ)$ coincide; namely, 
\[x\cdot y-y\cdot x=x\succ y+x\prec y-y\succ x-y\prec x=x\circ y-y\circ x\,.\]
If $f:D\to D'$ is a morphism of dendriform algebras, then it is also a morphism with
respect to any of the other three structures on $D$. The situation may be summarized  by
means of the following commutative diagram of categories
\[\xymatrix{ {\dendri}\ar[d]\ar[r] & {\prelie}\ar[d] \\
 {\assoc}\ar[r] & {\lie} } \]

Recall the notion of quasitriangular $\infbis$ from Section~\ref{S:infmod}.
In this section we show that there is a commutative diagram as follows:
\begin{equation}\label{E:quasi-diagram}
\xymatrix@1{ {\mbox{Quasitriangular
$\infbis$}}\ar[d]\ar[r] & {\infbis}\ar[d]\\ {\dendri}\ar[r] & {\prelie}  } 
\end{equation}
In this diagram,  the right vertical arrow is the functor
constructed in Section~\ref{S:prelie}, the bottom horizontal arrow is the
construction just discussed~\eqref{E:dendprelie}
 and the top horizontal arrow is simply the inclusion. It remains to discuss the
construction of a dendriform algebra from a quasitriangular $\infbi$, and to verify the
commutativity  of the diagram.

This construction is best understood from the point of view of {\em Baxter operators}.

\bd\label{D:baxter} Let $A$ be an associative algebra. A Baxter operator is a map
$\beta:A\to A$ that satisfies the condition
\begin{equation}\label{E:baxter}
\beta(x)\beta(y)=\beta\Bigl(x\beta(y)+\beta(x) y\Bigr)\,.
\end{equation}
\ed
Baxter operators arose in probability theory~\cite{Bax} and were a subject of
interest to Gian-Carlo Rota~\cite{R1,R2}.

We start by recalling a basic result from~\cite{poisson}, which provides us with the
 examples of Baxter operators that are most relevant for our present purposes. 
\bp\label{P:ayb-baxter}. Let   $r=\sum_i u_i\ten v_i$ be a solution of the associative
Yang-Baxter equation~\eqref{E:ayb} in an associative algebra $A$. Then the map $\beta:A\to
A$ defined by
\[\beta(x)=\sum_i   u_ixv_i\]
is a Baxter operator. 
\ep
\bpf Replacing the tensor
symbols in the associative Yang-Baxter equation~\eqref{E:ayb} by $x$ and $y$ one obtains
precisely~\eqref{E:baxter}.
\epf
\begin{rem}\label{R:baxter} The associative Yang-Baxter equation is analogous to the {\em
classical Yang-Baxter equation}, which is named after C. N. Yang and R. J. Baxter. Baxter
operators, on the other hand, are named after Glen Baxter.
\end{rem}

The following result provides the second step in the construction of dendriform algebras
from quasitriangular $\infbis$. 
\bp\label{P:baxter-dendri}
Let $A$ be an associative algebra and $\beta:A\to A$ a Baxter operator.
 Define  new operations on $A$ by
\[x\succ y=\beta(x) y \text{ \ and \ } x\prec y=x\beta(y)\ .\]
Then $(A,\succ ,\prec )$ is a dendriform algebra.
\ep
\bpf We verify the last axiom in~\eqref{E:dendri}; the others are similar. We have
\begin{align*}
x\succ (y\succ z) &=\beta(x)(y\succ z)=\beta(x)\beta(y)z\\
&\equal{\eqref{E:baxter}}\beta\Bigl(x\beta(y)+\beta(x) y\Bigr)z
=\beta(x\prec y+x\succ y)z\\ &=\beta(x\prec y)z+\beta(x\succ y)z
=(x\prec y)\succ z+(x\succ y)\succ z\,.
\end{align*}
\epf
Extensions of the above result appear in~\cite[Propositions 5.1 and 5.2]{poisson}.

Finally, we have the desired construction of dendriform algebras from quasitriangular
$\infbis$. 
\bt\label{T:quasi-dendri} Let $(A,r)$ be a quasitriangular $\infbi$, $r=\sum_i u_i\ten
v_i$. Define new operations on $A$ by 
\begin{equation}\label{E:quasi-dendri}
x\succ y=\sum_i u_ixv_iy \text{ \ and \ }x\prec y=\sum_i xu_iyv_i\,.
\end{equation}
Then, $(A,\succ,\prec)$ is a dendriform algebra.
\et
\bpf Combine Propositions~\ref{P:ayb-baxter} and~\ref{P:baxter-dendri}.
\epf

A morphism between quasitriangular $\infbis$ $(A,r)$ and $(A',r')$ is a morphism of
algebras $f:A\to A'$ such that $(f\ten f)(r)=r'$. Clearly, such a map $f$ preserves the
dendriform structures on $A$ and $A'$. Thus, we have constructed a
functor from quasitriangular $\infbis$ to dendriform algebras. 

We briefly discuss the functoriality of the construction with respect to derivations.

A derivation of a quasitriangular $\infbi$ $(A,r)$ is a map $D:A\to A$ that is a derivation
of the associative algebra $A$ such that
\[(D\ten\id+\id\ten D)(r)=0\,.\]
This implies that $D$ is a biderivation of the $\infbi$ associated to $(A,r)$, in
 the sense of~\eqref{E:bider}. In fact, it is easy to see that $D$ is a biderivation if and
only 
$(D\ten\id+\id\ten D)(r)$ is an invariant element in the \odash{A}{bimodule} $A\ten A$.
These two conditions are analogous to the ones encountered in the definition of
quasitriangular and coboundary $\infbis$~\cite[Section 1]{analog}. The stronger condition
guarantees the following:

\bp\label{P:quasi-der} Let $D$ be a derivation of a quasitriangular $\infbi$ $(A,r)$.
Then $D$ is also a derivation of the associated dendriform algebra, i.e.,
\[D(a\succ b)=a\succ D(b)+D(a)\succ b \text{ \ and \ }D(a\prec b)=a\prec D(b)+D(a)\prec
b\,.\]
\ep
\bpf Similar to the other proofs in this section.
\epf

It remains to verify the
commutativity of diagram~\eqref{E:quasi-diagram}. Starting from $(A,r)$ and going 
clockwise, we pass through the $\infbi$ with comultiplication $\Delta(b)=\sum_i u_i\ten
v_ib-\sum_i bu_i\ten v_i$, according to~\eqref{E:principal} (see also
Remark~\ref{R:signs}). The associated pre-Lie algebra structure is, by~\eqref{E:inflie},
$a\circ b=\sum_i u_iav_ib-\sum_i bu_iav_i$.
According to~\eqref{E:quasi-dendri}, this expression is equal to
$a\succ b-b\prec a$, which by~\eqref{E:dendprelie} is the pre-Lie algebra structure
obtained by going counterclockwise around the diagram.

\bex
\be
\item Let $A$ be an associative unital algebra and $b\in A$ an element such that
$b^2=0$. Then, $r:=1\ten b$ is a solution of the associative Yang-Baxter
equation~\eqref{E:ayb}~\cite[Example 5.4.1]{cor}. 
The corresponding dendriform structure on $A$ is simply
\[x\succ y=xyb \text{ \ and \ }x\prec y=xby\,.\]
This structure is well defined even if $A$ does not have a unit.
\smallskip
\item The element $r:=\mat{1 & 0\\0 & 0}\ten\mat{0 & 1\\0 & 0}-
\mat{0 & 1\\0 & 0}\ten\mat{1 & 0\\0 & 0}$ is a solution of~\eqref{E:ayb} in the
algebra
\smallskip
of matrices $M_2(k)$~\cite[Example 5.4.5.e]{cor} and~\cite[Examples
2.3.1 and 2.8]{analog}. The corresponding dendriform structure on $M_2(k)$ is 
\[\mat{a & b\\c & d}\succ\mat{x & y\\z & w}=\mat{az-cx & aw-cy\\ 0 & 0}\text{ \ and \ }
\mat{a & b\\c & d}\prec\mat{x & y\\z & w}=\mat{-az & ax\\-cz & cx}\,.\]

\item The most important example is provided by Drinfeld's double, which is a
quasitriangular $\infbi$ canonically associated to an arbitrary finite dimensional
$\infbi$. This will occupy the rest of the section. 
\ee
\eex

Let $A$ be a finite dimensional $\infbi$. Recall the definition of the double  $D(A)$
from Section~\ref{S:infmod}. Let $\{e_i\}$ be a linear basis of $A$ and $\{f_i\}$ the dual
basis of $A^*$. Let $r\in D(A)\ten D(A)$ be the element
\[r=\sum_i e_i\ten f_i\in A\ten A^*\subseteq D(A)\ten D(A)\,.\]
According to~\cite[Theorem 7.3]{cor}, $(D(A),r)$ is a quasitriangular $\infbi$ (see
Remark~\ref{R:signs}). By Theorem~\ref{T:quasi-dendri}, there is a dendriform algebra
structure on the space $D(A)=(A\ten A^*)\dirsum A\dirsum A^*$. 

In order to make this structure explicit, we introduce some notation. We identify
 $A\ten A^*$ with $\End(A)$ via 
\[(a\tie f)(b)=f(b)a\,.\]
For each $a\in A$ and $f\in A^*$, define linear endomorphisms of $A$ by
\begin{align*}
L_a(x)&=ax & & L_f(x)=f(x_2)x_1\\
R_a(x)&=xa & & R_f(x)=f(x_1)x_2\\
P_a(x)&=x_1ax_2 & & P_f(x)=f(x_2)x_1x_3\,.
\end{align*}

The composition of linear maps $\phi:U\to V$ and $\psi:V\to W$ is denoted by
$\psi\phi$, or $\psi(\phi)$, if the expression for $\phi$ is complicated. This should not
be confused with the evaluation of an endomorphism $T$ on an element $a\in A$, denoted by
$T(a)$. The convolution of linear endomorphisms $T$ and $S$ of $A$ is 
\[T\ast S=\mu(T\ten S)\Delta\,.\]

\bt\label{T:double-dendri} Let $A$ be an arbitrary $\infbi$. There is a dendriform
structure on the space $\End(A)\dirsum A\dirsum A^*$, given explicitly as follows. For
$a$, $b\in A$, $f$, $g\in A^*$ and $T$, $S\in\End(A)$,
\begin{align*}
a\succ b&=P_a(b)+R_aL_b & &a\prec b=L_aR_b\\
f\succ a&=P_f(a)+L_fL_a & &a\prec f=L_aL_f\\
f\prec a&=fP_a+R_fR_a & & a\succ f=R_aR_f\\
f\prec g&=fP_g+R_fL_g & & f\succ g=L_fR_g
\end{align*}
\begin{align*}
a\succ T&=P_aT+R_a(T\ast\id) & &a\prec T=L_a(\id\ast T)\\
f\succ T&=P_fT+L_f(T\ast\id) & &f\prec T=f(\id\ast T\ast\id)+R_f(\id\ast T)\\
T\prec a&=TP_a+(T\ast\id)R_a & & T\succ a=(\id\ast T\ast\id)(a)+(\id\ast T)L_a\\
T\prec f&=TP_f+(T\ast\id)L_f & & T\succ f=(\id\ast T)R_f
\end{align*}
\begin{align*}
T\succ S &=(\id\ast T\ast id)S +(\id\ast T)(S\ast \id)\\
T\prec S &=T(\id\ast S\ast\id)+(T\ast \id)(\id\ast S)
\end{align*}
\et
\bpf Assume that $A$ is finite dimensional, so $D(A)$ and $r=\sum e_i\ten f_i$ are well
defined, and hence there is a dendriform structure on the space $D(A)$.  For the details of the
infinite dimensional case see Remark~\ref{R:inf-proof}. 

We provide the derivations of the first and last formulas, the others are similar. 
We make use of the $\infbi$ structure of $D(A)$ as described in Proposition~\ref{P:double}.

For the first formula we have
\begin{align*}
a\succ b &=\sum_i e_iaf_ib=\sum_i e_ia(f_i\to b)+\sum_i e_ia(f_i\gets b)\\
&\equal{\eqref{E:actionone}}\sum_i e_iaf_i(b_1)b_2+\sum_i e_ia\tie (f_i\gets b)\\
&=b_1ab_2+\sum_i e_ia\tie (f_i\gets b)\,.
\end{align*}
Now, for any $x\in A$,
\[\sum_i (f_i\gets b)(x)e_ia\equal{\eqref{E:actiontwo}}\sum_i f_i(bx)e_ia=bxa=R_aL_b(x)\,.\]
Thus, $a\succ b=P_a(b)+R_aL_b$, as claimed.

For the last formula, let $T=a\tie f$ and $S=b\tie g$. We have
\begin{align*}
T\prec S&=(a\tie f)\prec (b\tie g)=\sum_i (a\tie f)e_i(b\tie g)f_i\\
&=\sum_i a(f\to e_ib)\tie gf_i+\sum_i a\tie(f\gets e_ib)gf_i\,.
\end{align*}
Hence, for any $x\in A$,
\begin{align*}
(T\prec S)(x)&\equal{\eqref{E:actiontwo}}
\sum_i f_i(x_1)g(x_2)a(f\to e_ib)+\sum_i f_i(x_1)g(x_2)(f\gets e_ib)(x_3)a\\
&=g(x_2)a(f\to x_1b)+g(x_2)(f\gets x_1b)(x_3)a\\
&\equal{\eqref{E:actiontwo}}a\bigl(f\to x_1S(x_2)\bigr)+g(x_2)f(x_1bx_3)a\\
&\equal{\eqref{E:actionone}}f\Bigl(\bigl(x_1S(x_2)\bigr)_1\Bigr)a\bigl(x_1S(x_2)\bigr)_2
+f\bigl(x_1S(x_2)x_3\bigr)a\\
&=(T\ast\id)(\id\ast S)(x)+T(\id\ast S\ast\id)(x)\,.
\end{align*}
Thus, $T\prec S=(T\ast\id)(\id\ast S)+T(\id\ast S\ast\id)$, as claimed.
\epf
\begin{rem}\label{R:inf-proof}
The $\infbi$ structure on $D(A)$ and the element $r\in D(A)\ten D(A)$ are well defined
only if $A$ is finite dimensional. However, all formulas in Theorem~\ref{T:double-dendri}
make sense and the theorem holds even if $A$ is infinite dimensional. This may
be seen as follows. There is always an algebra structure on the space $\End(A)\dirsum A\dirsum
A^*$, extending that of $D(A)$. Moreover, there is always a Baxter operator on this
algebra,  well defined by
\[\beta(a)=R_a,\ \ \beta(f)=L_f \text{ \ and \ }\beta(T)=\id\ast T\,.\]
It is easy to see that this coincides with the operator corresponding to
$r$, when $A$ is finite dimensional. In the general case, it  may be checked directly
that $\beta$ satisfies~\eqref{E:baxter}. The result then follows from
Proposition~\ref{P:baxter-dendri}.
\end{rem}
\begin{rem} In order to fully appreciate the symmetry in the previous formulas,
the following relations should be kept in mind:
\begin{align*}
T\ast R_a&=R_a(T\ast\id) & & & T\ast L_f&=(T\ast\id)L_f\\
L_a\ast T&=L_a(\id\ast T) & & & R_f\ast T&=(\id\ast T)R_f\\
L_aR_b&=R_bL_a & & & L_fR_g&=R_gL_f
\end{align*}
\end{rem}
\begin{rem} Consider the pre-Lie algebra structures on $A$ and $A'$ corresponding to their
$\infbi$ structures by means of Theorem~\ref{T:inflie}.
Since $A$ and $A'$ are $\infbisubs$ of $D(A)$ (Proposition~\ref{P:double}), the functoriality
of the construction implies that $A$ and $A'$ are pre-Lie subalgebras of $D(A)$, with respect
to the pre-Lie structure associated to the dendriform structure as in~\eqref{E:dendprelie}.
Let us verify this fact explicitly. The pre-Lie structure on $A$ is
\[a\circ b=a\succ b-b\prec a=P_a(b)+R_aL_b-L_bR_a=P_a(b)=b_1ab_2\,,\]
as expected. The pre-Lie structure on $A'$ is
\[f\circ g=f\succ g-g\prec f=L_fR_g-gP_f-R_gL_f=-gP_f\,.\]
Thus, 
\[(f\circ g)(a)=-g(f(a_2)a_1a_3)\equal{\eqref{E:coproddual}}f(a_2)g_2(a_1)g_1(a_3)\equal
{\eqref{E:proddual}} (g_1fg_2)(a)\,,\] 
also as expected (the $\infbi$ structure on $A'$ was
described in Section~\ref{S:infmod}).
\end{rem}

When $A$ is a quasitriangular $\infbi$, there are dendriform algebra structures both on $A$ and
$\End(A)\dirsum A\dirsum A^*$, by Theorems~\ref{T:quasi-dendri} and~\ref{T:double-dendri}. These
two structures are related by a canonical morphism of dendriform algebras.

\bp\label{P:pi-dendri}
Let $(A,r)$ be a quasitriangular $\infbi$, $r=\sum_i u_i\ten v_i$. Then, the map
\[\pi:\End(A)\dirsum A\dirsum A^*\to A,\ \ 
\pi(a)=a,\ \pi(f)=\sum_i f(u_i)v_i \text{ and } \pi(T)=\sum_i T(u_i)v_i\]
is a morphism of dendriform algebras.
\ep
\bpf  Assume that $A$ is finite dimensional.  According to~\cite[Proposition 7.5]{cor}, the
above formulas define a morphism of $\infbis$ $\pi:D(A)\to A$  (see Remark~\ref{R:signs}).
By the functoriality of the construction of dendriform algebras, $\pi$ is also a morphism
of dendriform algebras.

The general case may be obtained by showing that $\pi$ commutes with the Baxter operators on
$\End(A)\dirsum A\dirsum A^*$ and $A$. This follows from~\eqref{E:ayb}
and~\eqref{E:deltaonr}. 
\epf

The formulas in Theorem~\ref{T:double-dendri} show that $\End(A)$ is closed under
the dendriform operations. Together with Proposition~\ref{P:pi-dendri}, this gives the
following: 
\bc\label{C:end-dendri} Let $A$ be an arbitrary $\infbi$. Then there is a dendriform algebra
structure on the space $\End(A)$ of linear endomorphisms of $A$, defined by
\[T\succ S =(\id\ast T\ast id)S +(\id\ast T)(S\ast \id) \text{ \ and \ }
T\prec S =T(\id\ast S\ast\id)+(T\ast \id)(\id\ast S)\,.\]
Moreover, if $A$ is quasitriangular, with $r=\sum_i u_i\ten v_i$, then there is a morphism of
dendriform algebras $\End(A)\to A$ given by
\[T\mapsto \sum_i T(u_i)v_i\]
\ec
\bpf \epf

\begin{rem} There are in fact other, more primitive, dendriform structures on $\End(A)$
whenever $A$ is an $\infbi$. These will be studied in future work.
\end{rem}

\section{Infinitesimal Hopf bimodules, pre-Lie bimodules, dendriform
bimodules}\label{S:bimod}

In previous sections, we have shown how to construct a pre-Lie algebra from an $\infbi$
and a dendriform algebra from a quasitriangular $\infbi$. These constructions are
compatible, in the sense of~\eqref{E:quasi-diagram}. In this section we extend these
constructions to the corresponding  categories of bimodules.

The first step is to define the appropiate notion of bimodules over $\infbis$. Recall
the notion of left infinitesimal Hopf modules from Definition~\ref{D:infmod}. Right
infinitesimal Hopf modules are defined similarly. We combine these two notions in the
following:
\bd\label{D:infbimod}
Let $(A,\mu,\Delta)$ be an $\infbi$. An infinitesimal Hopf bimodule (abbreviated
$\infbimod$) over
$A$ is a space $M$ endowed with maps
\[\lambda:A\ten M\to M,\ \Lambda:M\to A\ten M,\ \xi:M\ten A\to M \text{ and }
\Xi:M\to M\ten A\]
such that
\bi
\item[(a)] $(M,\lambda,\Lambda)$ is a left $\infmod$ over $(A,\mu,\Delta)$,
\item[(b)] $(M,\xi,\Xi)$ is a right $\infmod$ over $(A,\mu,\Delta)$,
\item[(c)] $(M,\lambda,\xi)$ is a bimodule over $(A,\mu)$,
\item[(d)] $(M,\Lambda,\Xi)$ is a bicomodule over $(A,\Delta)$, and
\item[(e)] the following diagrams commute:
\begin{equation*}\label{E:infbimod}
\xymatrix@1{ {A\ten M}\ar[r]^{\id\ten\Xi}\ar[d]_{\lambda} & {A\ten
M\ten A}\ar[d]^{\lambda\ten\id} & & {M\ten A}\ar[r]^{\Lambda\ten\id}\ar[d]_{\xi} & {A\ten M\ten
A}\ar[d]^{\id\ten\xi}\\ 
{M}\ar[r]_{\Xi} & {M\ten A} & & {M}\ar[r]_{\Lambda} & {A\ten M}  } 
\end{equation*}
\ei
\ed
\begin{exa} For any $\infbi$ $(A,\mu,\Delta)$, the space $M=A\ten A$ is an $\infbimod$ via
\[\lambda=\mu\ten\id,\ \Lambda=\Delta\ten\id,\ \xi=\id\ten\mu \text{ and
}\Xi=\id\ten\Delta\,.\]
Note that $A$ itself, with the canonical bimodule and bicomodule structures, is
not an $\infbimod$.
\end{exa}

We will often use the following notation, for an $\infbimod$ $(M,\lambda,\Lambda,\xi,\Xi)$:
\begin{equation}\label{E:bimodnotation}
\lambda(a\ten m)=am,\ \ \xi(m\ten a)=ma,\ \
\Lambda(m)=m_{-1}\ten m_0 \text{ \ and \ } \Xi(m)=m_0\ten m_1\,.
\end{equation}
As is well known, this notation efficiently encodes the bicomodule axioms. For instance,
$m_{-2}\ten m_{-1}\ten m_0$ stands for
$(\Delta\ten\id)\Lambda(m)=(\id\ten\Lambda)\Lambda(m)$, and 
$m_{-1}\ten m_{0}\ten m_1$ for $(\id\ten\Xi)\Lambda(m)=(\Lambda\ten\id)\Xi(m)$.

 Just as left $\infmods$ over $A$ are left modules over $D(A)$
(Theorem~\ref{T:moddouble}),  $\odash{\epsilon}{Hopf}$ bimodules over $A$ are bimodules
over $D(A)$.

\bp\label{P:bimoddouble}
Let $A$ be a finite dimensional $\infbi$ and $M$ a space.
 If $(M,\lambda,\Lambda,\xi,\Xi)$ is an $\infbimod$ over $A$ as in~\eqref{E:bimodnotation},
then $M$ is a bimodule over $D(A)$ via
\begin{gather}
a\cdot m= am,\ \ f\cdot m=f(m_{-1})m_0,\ \ (a\tie
f)\cdot m=f(m_{-1})am_0 \label{E:bimoddouble1}\
\intertext{and}
m\cdot a= ma,\ \  m\cdot f=f(m_{1})m_0,\ \ m\cdot (a\tie
f)=f(m_{1})m_0a\,.\label{E:bimoddouble2}
\end{gather}
 Conversely, if $M$ is a bimodule over $D(A)$ then $M$ is an $\infbimod$ over $A$
and the structures are related as above.
\ep
\bpf Suppose $M$ is an $\infbimod$ over $A$. Then $(M,\lambda,\Lambda)$ is a left $\infmod$
and $(M,\xi,\Xi)$ is a right $\infmod$ over $A$. By Theorem~\ref{T:moddouble},
$M$ is a left and right module over $D(A)$ by means of~\eqref{E:bimoddouble1}
and~\eqref{E:bimoddouble2}. It only remains to check that these structures commute.

By assumption (c) (resp. (d)) in Definition~\ref{D:infbimod}, the left action of
$A$ (resp. $A'$) commutes with the right action of $A$ (resp. $A'$).  Similarly,
by assumption (e), the left action of $A$ (resp. $A'$) commutes with the
right action of $A'$ (resp. $A$). Since, by Proposition~\ref{P:double}, $D(A)$ is generated
as an algebra by $A\dirsum A'$, the previous facts guarantee that the left and right
actions of $D(A)$ on $M$ commute.

The converse is similar.
\epf

Next, we will relate $\infbimods$ over an $\infbi$ $A$ to bimodules over the associated
pre-Lie algebra, and similarly for the case of a quasitriangular $\infbi$ and the
associated dendriform algebra. For this purpose, we recall the  definition
of bimodules over these types of algebras. There is a general notion from the
theory of operads that dictates the bimodule axioms in each case~\cite{FM}.
In the cases of present interest, it turns out that the bimodule axioms are obtained
from the axioms for the corresponding type of algebras 
 by the following simple procedure. Each  axiom for the
given type of algebras yields three bimodule axioms, obtained by choosing one of
the variables $x$, $y$ or $z$ and replacing it by a variable $m$ from the bimodule (this
may yield repeated axioms). This leads to the following definitions.

\bd\label{D:prelie-bimod}Let $(P,\circ)$ be a (left) pre-Lie algebra. A
$\odash{P}{bimodule}$ is a space $M$ endowed with maps $P\ten M\to M$, $x\ten m\mapsto
x\circ m$ and $M\ten P\to M$, $m\ten x\mapsto m\circ x$, such that
\begin{gather}
x\circ(y\circ m)-(x\circ y)\circ m=y\circ(x\circ m)-(y\circ x)\circ
m\,,\label{E:prelie-bimod1}\\
x\circ(m\circ z)-(x\circ m)\circ z=m\circ(x\circ z)-(m\circ x)\circ
z\,.\label{E:prelie-bimod2}
\end{gather}
\ed
In Section~\ref{S:prelie} we encountered left \odash{P}{modules}~\eqref{E:preliemod}.
Note that any such can be turned into a $\odash{P}{bimodule}$ by choosing the trivial right
action $m\circ x\equiv 0$.

\bd\label{D:dendri-bimod} Let $(D,\succ,\prec)$ be a dendriform algebra. A
\odash{D}{bimodule}  is a vector space $M$ together with four maps
\begin{align*}
D\ten M\to M & & D\ten M\to M & &  M\ten D\to M & & M\ten D\to M\\
x\ten m\mapsto x\succ m & & x\ten m\mapsto x\prec m & & m\ten x\mapsto m\succ x & &
m\ten x\mapsto m\prec x 
\end{align*}
such that
\begin{multline*}
(x\prec y)\prec m = x\prec (y\prec m)+x\prec (y\succ m)\,,\\  
x\succ (y\prec m) = (x\succ y)\prec m\,,\\ 
x\succ (y\succ m) =(x\prec y)\succ m+(x\succ y)\succ m\,,
\end{multline*}
\begin{multline*}
(x\prec m)\prec z = x\prec (m\prec z)+x\prec (m\succ z)\,,\\ 
x\succ (m\prec z) = (x\succ m)\prec z\,,\\
x\succ (m\succ z) =(x\prec m)\succ z+(x\succ m)\succ z\,,
\end{multline*}
\begin{multline*}
(m\prec y)\prec z = m\prec (y\prec z)+m\prec (y\succ z)\,,\\ 
m\succ (y\prec z) = (m\succ y)\prec z\,,\\
m\succ (y\succ z) =(m\prec y)\succ z+(m\succ y)\succ z\,.
\end{multline*}
\ed

It is easy to verify that if $M$ is dendriform bimodule over $D$ then it is also a 
pre-Lie bimodule over the associated pre-Lie algebra~\eqref{E:dendprelie}  by
means of
\begin{equation}\label{E:dendprelie-bimod}
x\circ m=x\succ m-m\prec x \text{ \ and \ }m\circ
x=m\succ x-x\prec m\,.
\end{equation}

Next we show that the construction of pre-Lie algebras from $\infbis$ can be
extended to bimodules. 
\bp\label{P:preliebimod} Let $A$ be an $\infbi$ and $(M,\lambda,\Lambda,\xi,\Xi)$ an
$\infbimod$ over $A$. Then $M$ is a  bimodule over the pre-Lie algebra
$(A,\circ)$ of Theorem~\ref{T:inflie} via
\begin{equation}\label{E:preliebimod}
a\circ m= m_{-1}am_0+m_0am_1 \text{ \ and \ } m\circ
a=a_1ma_2\,.
\end{equation}
\ep
\bpf Consider the first  axiom in Definition~\ref{D:prelie-bimod}. We have $b\circ
m=m_{-1}bm_0+m_0bm_1$. Using the axioms in Definition~\ref{D:infbimod} we calculate
\begin{gather*}
\Lambda(b\circ m)=m_{-2}\ten m_{-1}bm_0+m_{-1}b_1\ten b_2m_0+m_{-2}bm_{-1}\ten m_0+
m_{-1}\ten m_0b m_1\,,\\
\Xi(b\circ m)=m_{-1}bm_0\ten m_{1}+m_{0}b_1\ten b_2m_1+m_{0}bm_{1}\ten m_2+
m_{0}\ten m_1b m_2\,.\\
\end{gather*}
Hence,
\begin{align*}
a\circ(b\circ m) &=m_{-2}a m_{-1}bm_0+m_{-1}b_1a b_2m_0+m_{-2}bm_{-1}a m_0+
m_{-1}a m_0b m_1\\
&+m_{-1}bm_0a m_{1}+m_{0}b_1a b_2m_1+m_{0}bm_{1}a m_2+m_{0}a m_1b m_2\,.
\end{align*}
On the other hand, $a\circ b=b_1ab_2$, so
\[(a\circ b)\circ m=m_{-1}b_1ab_2m_0+m_0b_1ab_2m_1\,.\]
Therefore,
\begin{align*}
a\circ(b\circ m)-(a\circ b)\circ m&=m_{-2}a m_{-1}bm_0+m_{-2}bm_{-1}a m_0+
m_{-1}a m_0b m_1\\
&+m_{-1}bm_0a m_{1}+m_{0}bm_{1}a m_2+m_{0}a m_1b m_2\,.
\end{align*}
Since this expression is symmetric under $a\leftrightarrow b$,
Axiom~\eqref{E:prelie-bimod1} holds.

Consider now the second axiom.  We have $m\circ b=b_1mb_2$. Using the axioms in
Definition~\ref{D:infbimod} we calculate
$\Lambda(m\circ b)=b_1m_{-1}\ten m_0b_2+b_1\ten b_2mb_3$ and
$\Xi(m\circ b)=b_1m_0\ten m_1b_2+b_1mb_2\ten b_3$.
It follows that
\[a\circ(m\circ b) =b_1m_{-1}a m_{0}b_2+b_1a b_2mb_3+b_1m_{0}a m_1b_2+
b_1mb_2a b_3\,.\]
Since 
\[(a\circ m)\circ b=b_1m_{-1}am_0b_2+b_1m_0am_1b_2\,.\]
we deduce that
\begin{equation}\tag{$*$}
a\circ(m\circ b)-(a\circ m)\circ b=b_1a b_2mb_3+b_1mb_2a b_3\,.
\end{equation}
On the other hand, since
$\Delta(a\circ b)=\Delta(b_1ab_2)=b_1\ten b_2ab_3+b_1a_1\ten a_2b_2+b_1ab_2\ten b_3$,
we have that
\[m\circ(a\circ b)=b_1m b_2ab_3+b_1a_1ma_2b_2+b_1ab_2m b_3\,,\]
and since
\[(m\circ a)\circ b=b_1a_1ma_2b_2\,,\]
we obtain that,
\begin{equation}\tag{$**$}
m\circ(a\circ b)-(m\circ a)\circ b=b_1m b_2ab_3+b_1ab_2m b_3\,.
\end{equation}
Comparing $(*)$ with $(**)$ with we see that Axiom~\eqref{E:prelie-bimod2} holds as well.
\epf

\begin{rem} Proposition~\ref{P:preliemod} may be seen as the particular case of
Proposition~\ref{P:preliebimod} when the right module and comodule structures on $M$ are
trivial (i.e., zero).
\end{rem}

We set now to extend the construction of dendriform algebras from quasitriangular
$\infbis$ to the categories of bimodules. As in Section~\ref{S:dendri}, it is convenient to
consider the more general context of Baxter operators.

\bd\label{D:baxter-bimod} Let $A$ be an associative algebra and $\beta_A:A\to A$ a
Baxter operator~\eqref{E:baxter}. A Baxter operator on a \odash{A}{bimodule} $M$
(relative to $\beta_A$) is a map $\beta_M:M\to M$ such that
\begin{gather}\label{E:baxter-bimod1}
\beta_A(a)\beta_M(m)=\beta_M\Bigl(a\beta_M(m)+\beta_A(a) m\Bigr)\,,\\
\beta_M(m)\beta_A(a)=\beta_M\Bigl(m\beta_A(a)+\beta_M(m)
a\Bigr)\,.\label{E:baxter-bimod2}
\end{gather}
\ed

\bp\label{P:baxter-dendri-bimod}
Let $A$ be an associative algebra, $\beta_A:A\to A$ a Baxter operator on $A$, $M$ an
\odash{A}{bimodule} and $\beta_M$ a Baxter operator on $M$.
 Define  new actions of $A$ on $M$ by
\[a\succ m=\beta_A(a)m\,,\  m\succ a=\beta_M(m)a\,,\ a\prec m=a\beta_M(m) \text{ \ and \ }
m\prec a=m\beta_A(a)\,.\] 
Equipped with actions, $M$ is a bimodule over the dendriform algebra of
Proposition~\ref{P:baxter-dendri}.
\ep
\bpf Similar to the proof of Proposition~\ref{P:baxter-dendri}.
\epf

\bp\label{P:ayb-baxter-bimod}. Let   $r=\sum_i u_i\ten v_i$ be a solution of the
associative Yang-Baxter equation~\eqref{E:ayb} in an associative algebra $A$. Let $M$
be an \odash{A}{bimodule}. Then the map
$\beta_M:M\to M$ defined by
\[\beta_M(m)=\sum_i   u_imv_i\]
is a Baxter operator on $M$, relative to the Baxter operator on $A$ of
Proposition~\ref{P:ayb-baxter}
\ep
\bpf Replacing the first tensor
symbol in the associative Yang-Baxter equation~\eqref{E:ayb} by $a$ and the second by $m$
one obtains~\eqref{E:baxter-bimod1}. Replacing them in the other order
yields~\eqref{E:baxter-bimod2}.
\epf

Finally, we have the desired construction of dendriform bimodules from
bimodules over quasitriangular $\infbis$. 
\bc\label{C:quasi-dendri-bimod} Let $(A,r)$ be a quasitriangular $\infbi$, $r=\sum_i
u_i\ten v_i$. Let $M$ be an arbitrary \odash{A}{bimodule}. Define new actions of $A$ on $M$
by 
\begin{equation}\label{E:quasi-dendri-bimod}
a\succ m=\sum_i u_iav_im\,,\  m\succ a=\sum_i u_imv_ia\,,\ a\prec m=\sum_i au_imv_i \text{
\ and \ } m\prec a=\sum_i mu_iav_i\,.
\end{equation}
Equipped with these actions, $M$ is a bimodule over the dendriform algebra of
Theorem~\ref{T:quasi-dendri}.
\ec
\bpf Combine Propositions~\ref{P:ayb-baxter-bimod} and~\ref{P:baxter-dendri-bimod}.
\epf

We have thus constructed, from an $\infbimod$ over an $\infbi$, a bimodule
over the associated pre-Lie algebra (Proposition~\ref{P:preliebimod}), and from a
bimodule over a quasitriangular $\infbi$, a bimodule over the associated dendriform
algebra (Corollary~\ref{C:quasi-dendri-bimod}). Also, a bimodule over a dendriform
algebra always yields a bimodule over the associated pre-Lie
algebra~\eqref{E:dendprelie-bimod}. In order to close the circle, it remains to construct
an $\infbimod$ from a bimodule over a quasitriangular $\infbi$.

\bp\label{P:bimodquasi} Let $(A,r)$ be a quasitriangular $\infbi$ and $M$ an
\odash{A}{bimodule}. Then $M$ becomes a  $\infbimod$ over $A$ via $\Lambda:M\to A\ten M$
and $\Xi:M\to M\ten A$ defined by
\begin{equation}\label{E:bimodquasi}
\Lambda(m)=\sum_i u_i\ten v_im \text{ \ and \ }\Xi(m)=-\sum_i mu_i\ten v_i\,.
\end{equation}
\ep
\bpf We already know, from Proposition~\ref{P:modquasi}, that $\Lambda$ turns $M$ into
a left $\infmod$ over $A$. Thus, axiom (a) in Definition~\ref{D:infbimod} holds.
Axiom (b) can be checked similarly: $\Xi$ is coassociative because of the fact that
$(\id\ten\Delta)(r)=-r_{13}r_{12}$, which holds according to~\cite[Proposition 5.5]{cor}
(one must take into account Remark~\ref{R:signs}). Note that the minus sign in the
definition of $\Xi$ is essential, to ensure both this and the right $\infmod$
axiom.

Axiom (c) holds by hypothesis. Axiom (d) holds because both $(\id\ten\Xi)\Lambda(m)$
and $(\Lambda\ten\id)\Xi(m)$ are equal to
\[-\sum_{i,j}u_i\ten v_imu_j\ten v_i\,.\]
Axiom (e) holds because both $\Xi\lambda(a\ten m)$ and $(\lambda\ten\id)(\id\ten\Xi)(a\ten
m)$ are equal to 
\[-\sum_i amu_i\ten v_i\,,\]
 and similarly for $\Lambda\xi$ and
$(id\ten\xi)(\Lambda\ten\id)$. This completes the proof.
\epf

Let $(A,r)$ be a quasitriangular $\infbi$. Consider the associated algebras as
in diagram ~\eqref{E:quasi-diagram}. In this section, we have constructed the corresponding
diagram at the level of bimodules:
\begin{equation}\label{E:quasi-diagram-bimod}
\xymatrix@1{ {\mbox{(Associative) bimodules}}\ar[d]\ar[r] & {\infbimods}\ar[d]\\
{\mbox{Dendriform bimodules}}\ar[r] & {\mbox{Pre-Lie bimodules}}  } 
\end{equation}
Each arrow is a functor, as morphisms are clearly preserved. The diagram is indeed
commutative. Going around clockwise, we pass through the $\infbimod$ with coactions
$\Lambda(m)=\sum_i u_i\ten v_im$ and $\Xi(m)=-\sum_i mu_i\ten v_i$, according
to~\eqref{E:bimodquasi}. The associated pre-Lie bimodule actions are,
by~\eqref{E:preliebimod},
\[a\circ m=\sum_i u_iav_im-\sum_i mu_iav_i \text{ \ and \ } 
m\circ a=\sum_i u_imv_ia-\sum_i au_imv_i\,.\]
According to~\eqref{E:quasi-dendri-bimod}, these expressions are
respectively equal to
$a\succ m-m\prec a$ and $m\succ a-a\prec m$, which by~\eqref{E:dendprelie-bimod} is the
pre-Lie bimodule structure obtained by going counterclockwise around the diagram.

\section{Brace algebras}\label{S:brace}
In this section we explain how one may associate a brace algebra to an arbitrary $\infbi$, in a
way that refines the pre-Lie algebra construction of Section~\ref{S:prelie} and that is
compatible with the dendriform algebra construction of Section~\ref{S:dendri}.

We provide the left version of the definition of brace algebras given in~\cite{C}.
\bd\label{D:brace} A (left) brace algebra is a space $B$ equipped with multilinear operations
$B^{n}\times B\to B$, $(x_1,\ldots,x_n,z)\mapsto \brac{x_1,\ldots,x_n;z}$, one for
each $n\geq 0$, such that
\[\brac{z}=z\]
and for any $n$, $m\geq 1$,
\begin{equation}
\begin{split}
\brac{x_1,\ldots,x_n;\brac{y_1,\ldots,y_m;z}}&= \label{E:brace}\\
\sum&\brac{X_0,\brac{X_1;y_1},X_2,\brac{X_3;y_2},X_4,\ldots,X_{2m-2},\brac{X_{2m-1};y_m},
X_{2m};z }\,,
\end{split}
\end{equation}
where the sum takes place over all partitions of the ordered set $\{x_1,\ldots,x_n\}$ into
(possibly empty) consecutive intervals $X_0\sqcup X_1\sqcup\cdots\sqcup X_{2m}$.
\ed

The case $n=m=1$ of Axiom~\ref{E:brace} says
\begin{equation}\label{E:braceone}
\brac{x;\brac{y;z}}=\brac{x,y;z}+
\brac{\brac{x;y};z}+\brac{y,x;z}\,.
\end{equation}
 The three terms on the right hand side correspond
respectively to the partitions
$(\{x\},\emptyset,\emptyset)$, $(\emptyset,\{x\},\emptyset)$ and $(\emptyset,\emptyset,\{x\})$.

The operation $x\circ y:=\brac{x;y}$ endows $B$ with a pre-Lie algebra structure. In fact,
\eqref{E:braceone} shows that $x\circ(y\circ z)-(x\circ y)\circ z$ is symmetric under
$x\leftrightarrow y$, so Axiom~\eqref{E:prelie} holds. This defines a functor from brace
algebras to pre-Lie algebras. The construction of pre-Lie algebras from $\infbis$ in
Section~\ref{S:prelie} can be refined accordingly, as we explain next.

In order to describe this
refined construction, we must depart from our notational convention for coproducts and revert to
Sweedler's original notation. Thus, in this section, the coproducts of an element $b$ will be
denoted by
\[\Delta(b)=\sum_{(b)}b_{(1)}\ten b_{(2)}\]
 and the \odash{n}{th} iteration of the coproduct  by
\[\Delta^{(n)}(b)=\sum_{(b)}b_{(1)}\ten\cdots\ten b_{(n+1)} \,.\]

\bt\label{T:infbrace} Let $A$ be an $\infbi$. Define operations $A^n\times A\to A$ by
\[\brac{a_1,\ldots,a_n;b}=\sum_{(b)} b_{(1)}a_1b_{(2)}a_2\ldots b_{(n)}a_nb_{(n+1)}\,.\]
These operations turn $A$ into a brace algebra.
\et
\bpf The complete details of the proof will be provided elsewhere. The idea is simple:
each term on the right hand side of~\eqref{E:brace} corresponds to a term in the expansion
of
\[\Delta^{(n)}\Bigl(\sum_{(z)}z_{(1)}y_1z_{(2)}y_2\ldots z_{(m)}y_m z_{(m+1)}\Bigr)\]
obtained by successive applications of~\eqref{E:infbi}. For instance, when $n=2$ and $m=1$,
one has
\begin{align*}
\Delta^{(2)}\Bigl(\sum_{(z)}z_{(1)}yz_{(2)}\Bigr)&
=\sum_{(z)}z_{(1)}\ten z_{(2)}\ten z_{(3)}yz_{(4)}+z_{(1)}\ten z_{(2)}y
z_{(3)}\ten z_{(4)}+z_{(1)}y z_{(2)}\ten z_{(3)}\ten z_{(4)}\\
&+\sum_{(z),(y)}z_{(1)}y_{(1)}\ten y_{(2)}z_{(2)}\ten z_{(3)}
+z_{(1)}y_{(1)}\ten y_{(2)}\ten y_{(3)}z_{(2)}+z_{(1)}\ten z_{(2)}y_{(1)}\ten y_{(2)}z_{(3)}
\end{align*}
Therefore,
\begin{align*}
\brac{x_1,x_2;\brac{y;z}}&
=\sum_{(z)}z_{(1)}x_1 z_{(2)}x_2 z_{(3)}yz_{(4)}+z_{(1)}x_1 z_{(2)}yz_{(3)}x_2 z_{(4)}
+z_{(1)}y z_{(2)}x_1 z_{(3)}x_2 z_{(4)}\\
&+\sum_{(z),(y)}z_{(1)}y_{(1)}x_1y_{(2)}z_{(2)}x_2z_{(3)}
+z_{(1)}y_{(1)}x_1y_{(2)}x_2y_{(3)}z_{(2)}+z_{(1)}x_1z_{(2)}y_{(1)}x_2y_{(2)}z_{(3)}\\
&=\brac{x_1,x_2,y;z}+\brac{x_1,y,x_2;z}+\brac{y,x_1,x_2;z}\\
&+\brac{\brac{x_1;y},x_2;z}+\brac{\brac{x_1,x_2;y};z}+\brac{x_1,\brac{x_2;y};z}
\end{align*}
which is Axiom~\eqref{E:brace}.
\epf

By construction, the first brace operation on $A$ is simply
\[\brac{a;b}=\sum_{(b)}b_{(1)}ab_{(2)}\,,\]
which agrees with the pre-Lie operation~\eqref{E:inflie}. In this sense, the constructions
of Theorems~\ref{T:inflie} and~\ref{T:infbrace} are compatible.

\begin{exa} \label{Ex:brace}
Consider the $\infbi$ $k[\bfx,\bfx^{-1}]$ of divided differences
(Examples~\ref{Ex:prelie}). It is easy to see that for any $n\geq 0$ and $r\in\Z$,
\[\mu^{(n)}\Delta^{(n)}(\bfx^r)=\binom{r}{n}\bfx^{r-n}\,,\]
where it is understood, as usual, that
$\binom{r}{n}=0$ if $n>r\geq 0$ and $\binom{r}{n}=(-1)^n\binom{{-}r{+}n-1}{n}$ if $r<0$. It
follows that the brace algebra structure on $k[\bfx,\bfx^{-1}]$ is
\[\brac{\bfx^{p_1},\ldots,\bfx^{p_n};\bfx^r}=\binom{r}{n}\bfx^{r+p_1+\cdots+p_n-n}\,. \]
The brace axioms~\eqref{E:brace} boil down to a set of interesting identities
involving binomial coefficients.

Fr\'ed\'eric Chapoton made us aware of the fact that if one applies the general
construction of~\cite[Proposition 1]{GV95b} (dropping all signs)
to the associative operad, one obtains precisely the
brace subalgebra $k[\bfx]$ of our brace algebra  $k[\bfx,\bfx^{-1}]$.

This example may be generalized in another direction. Namely, if $A$ is a commutative
algebra and $D:A\to A$ a derivation, then one obtains a brace algebra structure on $A$ by
defining
\[\brac{x_1,\ldots,x_n;z}=x_1\cdots x_n\frac{D^n(z)}{n!} \]
(assuming $\ch(k)=0$). The example above corresponds to $A=k[\bfx,\bfx^{-1}]$,
$D=\frac{d}{d\bfx}$.
\end{exa}
\smallskip
Brace algebras sit between dendriform and pre-Lie algebras: Ronco has shown that one
can associate a brace algebra to a dendriform algebra, by means of certain
operations~\cite[Theorem 3.4]{R}. Our constructions of dendriform and brace algebras
from Theorems~\ref{T:quasi-dendri} and~\ref{T:infbrace} are compatible with this functor.

In summary, one obtains a commutative diagram
\[\xymatrix@1{ {\mbox{Quasitriangular
$\infbis$}}\ar[d]\ar[r] & {\infbis}\ar[d]\ar[rd]\\ {\dendri}\ar[r] & {\bracealg}\ar[r] &
{\prelie}  }\]
The details will be provided elsewhere.

\appendix

\section{Infinitesimal bialgebras as comonoid objects}\label{S:comonoid}

Ordinary bialgebras are bimonoid objects in the braided monoidal category of vector spaces,
where the monoidal structure is the usual tensor product $V\ten W$ and the braiding is the
trivial symmetry $x\ten y\mapsto y\ten x$. In this appendix, we construct a certain
monoidal category of algebras for which the comonoid objects are precisely $\infbis$.
Related notions of bimonoid objects are discussed as well.

For the basics on monoidal categories the reader is referred to~\cite[Chapter XI]{Kas}.
The monoidal categories we consider possess a unit object, and whenever we refer to monoid
objects these are assumed to be unital, even if not explicitly stated. Similarly, comonoid
objects are assumed to be counital.

We start by recalling the well known {\em circular tensor product} of vector spaces.

\bd \label{D:circular} The circular tensor product of two vector spaces $V$ and $W$ is
\[V\cirprod W=V\dirsum W\dirsum (V\ten W)\,.\]
We denote the elements of this space by triples $(v,w,x\ten y)$. The circular tensor
product of maps $f:V\to X$ and $g:W\to Y$ is 
\[(f\cirprod g)(v,w,x\ten y)=\bigl(f(v),g(w),(f\ten g)(x\ten y)\bigr)\,.\]
\ed

Both spaces $(U\cirprod V)\cirprod W$ and $U\cirprod (V\cirprod W)$ can be canonically
identified with 
\[U\dirsum V\dirsum W\dirsum (U\ten V)\dirsum (U\ten W)\dirsum (V\ten W)\dirsum (U\ten
V\ten W)\,.\]
This gives rise to a natural isomorphism $(U\cirprod V)\cirprod W\cong U\cirprod (V\cirprod
W)$ which satisfies the pentagon for associativity. 
This endows the category of vector spaces with a monoidal structure, for which the unit
object is the zero space. We denote this monoidal category by $(\Spa,\cirprod, 0)$.

Let $(\Spa,\ten,k)$ denote the usual monoidal category of vector spaces, where
the monoid objects are unital associative algebras and the comonoid objects are counital
coassociative  coalgebras.
There is an obvious monoidal functor $\alpha:(\Spa,\cirprod,0)\to(\Spa,\ten,k)$ defined by
\[V\mapsto V\dirsum k\,.\]
It is the so called {\em augmentation} functor. 

Monoids and comonoids in $(\Spa,\cirprod,0)$
are easy to describe: they are, respectively, non
unital algebras and non counital coalgebras (Proposition~\ref{P:mon-comon}, below).
Monoids and comonoids are preserved by monoidal functors. In the present situation this
simply says that a non unital algebra can be canonically augmented into a unital
algebra, and similarly for  coalgebras.
 
\bp\label{P:mon-comon}
A unital monoid object in $(\Spa,\cirprod,0)$ is precisely an associative algebra, not
necessarily unital. A counital comonoid object is precisely a coassociative coalgebra, not
necessarily counital.
\ep
\bpf Let $(A,\mu)$ be an associative algebra, $\mu(a\ten a')=aa'$. Define a map
$\tilde{\mu}:A\cirprod A\to A$ by
\begin{equation}\label{E:multcir}
(a,a',x\ten x')\mapsto a+a'+xx'\,.
\end{equation}
Let $u:0\to A$ be the unique map. Then, the diagrams
\[\xymatrix{
A\cirprod A\cirprod A\ar[r]^{\tilde{\mu}\cirprod\id}\ar[d]_{\id\cirprod\tilde{\mu}} &
A\cirprod A\ar[d]^{\tilde{\mu}}\\ A\cirprod A\ar[r]_{\tilde{\mu}} & A} \text{ \ and \ }
\xymatrix{
0\cirprod A\ar[r]^{u\ten\id}\ar[dr]_{\cong} & A\cirprod A\ar[d]|{\tilde{\mu}} & A\cirprod
0\ar[l]_{\id\cirprod u}\ar[dl]^{\cong}\\
 & A }\]
commute. Thus, $(A,\tilde{\mu},u)$ is a unital monoid in $(\Spa,\cirprod,0)$. 

Conversely, if $(A,\tilde{\mu},u)$ is a unital monoid in $(\Spa,\cirprod,0)$, then
$\tilde{\mu}$ must be of the form~\eqref{E:multcir} for an associative
multiplication $\mu$ on $A$, by the commutativity of the diagrams above.

The assertion for comonoids is similar. The comultiplication $\Delta:A\to A\ten A$ is
related to the comonoid structure $\tilde{\Delta}:A\to A\cirprod A$ by
\begin{equation}\label{E:deltacir}
\tilde{\Delta}(a)= (a,a,\Delta(a))
\end{equation}
and $\epsilon:A\to 0$ is the unique map.
\epf

\begin{rem} It is natural to wonder if there is a braiding on the monoidal category
$(\Spa,\cirprod,0)$ for which the bimonoid objects are precisely $\infbis$. 
 We know of two braidings on $(\Spa,\cirprod,0)$. The corresponding notions of bimonoid
objects are briefly discussed next. Neither yields $\infbis$.
\be
\item For any spaces $V$ and $W$, consider the map $\sigma_{V,W}:V\cirprod W\to W\cirprod V$
defined by
\[(v,w,x\ten y)\mapsto (w,v,y\ten x)\,.\]
This family of maps clearly satisfies the axioms for a braiding on the monoidal
category $(\Spa,\cirprod,0)$. Under the monoidal functor
$\alpha$, the braiding $\sigma$ corresponds to the usual braiding on $(\Spa,\ten,k)$ (the
trivial symmetry). For this reason, a bimonoid object in $(\Spa,\cirprod,0,\sigma)$
can be canonically augmented into an ordinary bialgebra.

It follows from  Proposition~\ref{P:mon-comon}) that a bimonoid object in
$(\Spa,\cirprod,0,\sigma)$ is a space $A$, equipped with an associative algebra structure
$A\ten A\to A$, $a\ten a'\mapsto aa'$, and a coassociative coalgebra structure
$\Delta:A\to A\ten A$, $a\mapsto a_1\ten a_2$, related by the axiom
\[\Delta(aa')=a\ten a'+a'\ten a+aa'_1\ten a'_2+a'_1\ten aa'_2+a_1a'\ten a_2+a_1\ten a_2a'+
a_1a'_1\ten a_2a'_2\,.\]
This is {\em not} the axiom which defines $\infbis$~\eqref{E:infbi}.

The axiom above is a translation of the fact that the map $A\to A\cirprod A$
 must be a morphism of monoids. We omit this calculation, but provide an explicit
description of the monoid structure on $A\cirprod A$.
More generally, we describe the tensor product of two monoids $A$ and $B$ in
$(\Spa,\cirprod,0,\sigma)$. 

According to Proposition~\ref{P:mon-comon}, the monoid structure on $A\cirprod B$
is uniquely determined by an associative multiplication on the space $A\cirprod B$.
We describe this multiplication, in terms of those of $A$ and $B$. It is
\[(a,b,x\ten y)\cdot(a',b',x'\ten y')=(aa',bb',a\ten b'+a'\ten b+ax'\ten y'+x'\ten by'+
xa'\ten y+x\ten yb'+xy\ten x'y')\,.\]
This is the result of composing
\[(A\cirprod B)\ten(A\cirprod B)\inc (A\cirprod B)\cirprod(A\cirprod
B)\map{\id\cirprod\sigma_{B,A}\cirprod\id} (A\cirprod A)\cirprod(B\cirprod
B)\map{\tilde{\mu}_A\cirprod\tilde{\mu}_B}A\cirprod B\,.\]

\item There exists a second braiding on the monoidal category $(\Spa,\cirprod,0)$,
for which bimonoid objects are somewhat closer to $\infbis$. It is the
family of maps 
$\beta_{V,W}:V\cirprod W\to W\cirprod V$ defined by
\[(v,w,x\ten y)\mapsto (w,v,0)\,.\]
Apart from the fact that $\beta$ is {\em not} an isomorphism, the
braiding axioms are satisfied by $\beta$. This allows us to construct a monoid structure
on the circular tensor product of two monoids in $(\Spa,\cirprod,0)$, and therefore to
speak of bimonoid objects in $(\Spa,\cirprod,0,\beta)$, as usual.

Since the monoidal functor $\alpha$ does not preserve this braiding, the
augmentation of a bimonoid in $(\Spa,\cirprod,0,\beta)$ is not an ordinary
bialgebra. Neither is it true that these bimonoid objects are $\infbis$.
In fact, a bimonoid object in
$(\Spa,\cirprod,0,\beta)$ is  a space $A$, equipped with an associative algebra structure
 and a coassociative coalgebra structure, as above, related by the axiom
\[\Delta(aa')=a\ten a'+aa'_1\ten a'_2+a_1\ten a_2a'\,.\]
Compare with Axiom~\eqref{E:infbi} for infinitesimal bialgebras.

This can be deduced  from the following description of the tensor product in
$(\Spa,\cirprod,0,\beta)$ of two monoid objects $A$ and $B$. This structure is
determined by the following (associative) multiplication on $A\cirprod B$:
\[(a,b,x\ten y)\cdot(a',b',x'\ten y')=(aa',bb',a\ten b'+ax'\ten y'+x\ten yb')\,.\]
 This is the result of composing
\[(A\cirprod B)\ten(A\cirprod B)\inc (A\cirprod B)\cirprod(A\cirprod
B)\map{\id\cirprod\beta_{B,A}\cirprod\id} (A\cirprod A)\cirprod(B\cirprod
B)\map{\tilde{\mu}_A\cirprod\tilde{\mu}_B}A\cirprod B\,.\]
\ee
\end{rem}

Let $\Alg$ denote the category of monoids in $(\Spa,\cirprod,0)$, that is,
associative algebras which are not necessarily unital.
We define a new monoidal structure on this category, independent of any braiding
on $(\Spa,\cirprod,0)$. We will show that $\infbis$ are precisely comonoid objects
in the resulting monoidal category. 

\bp\label{P:cirprod} Let $A$ and $B$ be associative algebras, not necessarily unital.
Then $A\cirprod B$ is an associative algebra via
\begin{equation}\label{E:cirprod}
(a,b,x\ten y)\cdot(a',b',x'\ten y')=(aa',bb',ax'\ten
y'+x\ten yb')\,.
\end{equation}
\ep
\bpf Consider the algebra $R=A\dirsum B$ and the \odash{R}{bimodule} $M=A\ten B$, with
\[(a,b)\cdot x\ten y=ax\ten y \text{ \ and \ } x\ten y\cdot (a,b)=x\ten yb\,.\]
The algebra $A\cirprod B$ is precisely the trivial extension $R\dirsum M$, where the
multiplication is
\[(r,m)\cdot(r',m')=(rr',rm'+mr')\,.\]
\epf

If $f:A\to A'$ and $g:B\to B'$ are morphisms of algebras, then so is $f\cirprod
g:A\cirprod B\to A'\cirprod B'$. In this way, $(\Alg,\cirprod,0)$ becomes a monoidal 
category. 

\bp\label{P:comonoidinalg}
A counital comonoid object in the monoidal category $(\Alg,\cirprod, 0)$ is
precisely an $\infbi$. 
\ep 
\bpf Let $(A,\mu,\Delta)$ be an $\infbi$. By Proposition~\ref{P:mon-comon},
$A$ may be seen as a monoid and comonoid in $(\Alg,\cirprod,0)$. It only remains to verify that
the comonoid structure $\tilde{\Delta}:A\to A\cirprod
A$, $\tilde{\Delta}(a)=(a,a,\Delta(a))$,  is a  morphism of algebras.
This is clear from~\eqref{E:cirprod} and~\eqref{E:infbi}. 

The converse is similar.
\epf 

The category of modules over an ordinary bialgebra $H$ is monoidal: the tensor
product of two \odash{H}{modules} acquires an \odash{H}{module} structure by restricting
the natural \odash{H\ten H}{module} structure via the comultiplication $\Delta:H\to H\ten
H$. There is no analogous construction for arbitrary $\infbis$. However, 
 it is possible to construct tensor products of certain modules over $\infbis$,
as discussed next.

\bp\label{P:cirmod} Let $A$ and $B$ be associative algebras. Let $M$ be a right
\odash{A}{module} and
$N$ a left \odash{B}{module}. Then $A\cirprod N$ is a left \odash{A\cirprod
B}{module} via
\[(a, b,x\ten y)\cdot(a',n,x'\ten v):=(aa',bn,ax'\ten v+x\ten yn)\,,\]
and $M\cirprod B$ is a right \odash{A\cirprod B}{module} via
\[(m,b',u\ten y')\cdot(a,b,x\ten y):=(ma,b'b,ub\ten y'+mx\ten y)\,.\]
\ep
\bpf Similar to the proof of Proposition \ref{P:cirprod}.
\epf

Let $A$ be an $\infbi$ and $N$ a left
\odash{A}{module}. It is possible to define a left \odash{A}{module} structure on
$A\cirprod N$, by restricting the structure of Proposition
\ref{P:cirmod} along the morphism of algebras $\tilde{\Delta}:A\to A\cirprod A$.
The resulting action of $A$ on $A\cirprod N$ is
\begin{equation}\label{E:ciraction}
a\cdot(a', n,x\ten v)=(aa',an+ax\ten v,a_1\ten a_2n)\,.
\end{equation}

\section{Counital infinitesimal bialgebras}\label{S:counital}

\bd\label{D:counital} An $\infbi$ $(A,\mu,\Delta)$ is said to be counital if 
the underlying coalgebra is counital, that is, if there exists a  map
$\eta:A\to k$ such that $(\id\ten\eta)\Delta=\id=(\eta\ten\id)\Delta$.
\ed
The map $\eta$ is necessarily unique and is called the counit of $A$. We use $\eta$
instead of the customary $\varepsilon$ to avoid confusion with the abbreviation for 
infinitesimal bialgebras.

Recall that if an $\infbi$ $A$ is both unital {\em and} counital then
$A=0$~\cite[Remark 2.2]{cor}. Nevertheless, many $\infbis$ arising in practice are either
unital or counital. In this appendix we study counital $\infbis$; all
constructions and results admit a dual version that applies to unital
$\infbis$.

We first show
that  counital $\infbis$ can be seen as comonoid objects in a certain monoidal category
of algebras. This construction is parallel to that for arbitrary $\infbis$ discussed
in Appendix~\ref{S:comonoid}. The two constructions are related  by means
of a pair of monoidal functors, but neither is more general than the other.

\bl\label{L:counital}
Let $A$ be a counital $\infbi$ with counit $\eta$. Then
\[\eta(aa')=0 \text{ for all $a$, $a'\in A$.}\]
\el
\bpf We show that  any coderivation $D:M\to C$ from a counital bicomodule $(M,s,t)$ to a
counital coalgebra
$(C,\Delta,\eta)$ maps to the kernel of $\eta$. The result follows by applying this remark to
the coderivation $\mu:A\ten A\to A$.

We have 
\begin{align*}
\eta D &=(\eta\ten\eta)\Delta D=(\eta\ten\eta)\bigl((\id_C\ten D) t+(D\ten\id_C)
s\bigr)\\
&= (\id_k\ten\eta D)(\eta\ten\id_M)t+ (\eta D\ten\id_k)(\id_M\ten\eta)s\\
&=\eta D+\eta D=2\cdot \eta D,
\end{align*}
whence $\eta D=0$.
\epf

This motivates the following definition.

\bd\label{D:augmented} Let $(A,\mu)$ be an algebra over $k$, not necessarily unital.
We say that it is {\em augmented} if there is given a map $\eta:A\to k$ such that
\[\eta(aa')=0 \text{ for all $a$, $a'\in A$.}\]
A morphism between augmented algebras $(A,\eta_A)$ and $(B,\eta_B)$ is a morphism of
algebras $f:A\to B$ such that $\eta_B f=\eta_A$.
\ed

\bp\label{P:tensorprod} Let $(A,\eta_A)$ and $(B,\eta_B)$ be augmented algebras.
Then $A\ten B$ is an associative algebra with multiplication
\begin{equation}
(a\ten b)\cdot(a'\ten b'):=\eta_B(b)aa'\ten b'+\eta_A(a')a\ten bb'\,.
\end{equation} 
Moreover, $A\ten B$ is augmented by
\[\eta_{A\ten B}(a\ten b):=\eta_A(a)\eta_B(b)\,.\]
\ep
\bpf The first assertion is Lemma 3.5.b in~\cite{cor} and the second is 
straightforward.
\epf

We denote the resulting augmented algebra by $A\teninf B$. This operation defines
a monoidal structure on the category of augmented algebras over $k$.
The unit object is the base field $k$ equipped with the zero multiplication and
the identity augmentation. We denote this monoidal category by $(\AAlg,\teninf, k)$. 

\bp\label{P:comonoid}
A counital comonoid object in the monoidal category $(\AAlg,\teninf, k)$ is
precisely a counital $\infbi$. 
\ep 
\bpf Start from a counital $\infbi$ $(A,\mu,\Delta,\eta)$. By Lemma \ref{L:counital},
$(A,\eta)$ is an augmented algebra. Moreover, by Lemma 3.6.b in~\cite{cor}, $\Delta:A\to
A\teninf A$ is a  morphism of algebras,  and it preserves the augmentations
by counitality. Clearly, $\eta:A\to k$ is also a  morphism of augmented algebras. 
Thus, $(A,\mu,\Delta,\eta)$ is a counital comonoid in $(\Alg,\teninf, k)$. 

Conversely, let $A$ be a counital comonoid in $(\AAlg,\teninf, k)$. First of all,
the counit $A\to k$ must preserve the augmentations of $A$ and $k$, 
 so it must coincide with the augmentation of $A$. The comultiplication must be
a morphism of algebras $A\to A\teninf A$. This implies Axiom~\eqref{E:infbi},
by definition of the algebra structure on $A\teninf A$ and counitality.
Thus $A$ is a counital $\infbi$.
\epf 

\begin{rem} An augmented algebra may be seen as a monoid in a certain
monoidal category of ``augmented vector spaces''. However, the monoidal structure
on
$(\AAlg,\teninf, k)$ does not come from a braiding on the larger category of augmented
vector spaces. For this reason, we cannot view counital $\infbis$ as bimonoid objects.
The situation parallels that encountered in Appendix~\ref{S:comonoid} for arbitrary
$\infbis$. In fact, there is pair of monoidal functors relating the two situations,
as we discuss next.
\end{rem}

\begin{rem}
Given a non unital algebra $A$, let $A^+:=A\dirsum k$, with algebra structure
\begin{equation}\label{E:multplus}
(a,x)\cdot(b,y)=(ab,0)\,.
\end{equation}
Note that $A^+$ is not the usual augmentation of $A$; in fact, $A^+$ is non unital.
Define $\eta(a,x)=x$. Then $\eta((a,x)\cdot(b,y))=0$, so $(A^+,\eta)$ is an
augmented algebra in the sense of Definition~\ref{D:augmented}.
Moreover, it is easy to see that there is a natural isomorphism of
augmented algebras
\[(A\cirprod B)^+\cong A^+\teninf B^+\,.\]
The application $A\mapsto A^+$ is thus a monoidal functor 
\[(\Alg,\cirprod, 0)\to(\AAlg,\teninf, k)\,.\]
The fact that comonoid objects are preserved by this monoidal functor simply says
that any $\infbi$ $A$ can be made into a counital $\infbi$ $A^+$, by extending the
comultiplication via $\Delta(1)=1\ten 1$ and the multiplication as in~\eqref{E:multplus}.

In the other direction, consider the forgetful functor
\[(\AAlg,\teninf, k)\to(\Alg,\cirprod, 0),\ \ (A,\eta)\mapsto A\,.\]
It is easy to see that the map
\begin{equation}\label{E:infcir}
A\teninf B\to A\cirprod B,\ \ a\ten b\mapsto
\eta_B(b)a+\eta_A(a)b+a\ten b
\end{equation}
 is a (natural) morphism of algebras. It follows that the
forgetful functor is {\em lax monoidal}. The fact that comonoid objects are preserved by
this type of functors
 simply says in this case that any counital $\infbi$ is in particular an $\infbi$.

Neither functor between these two categories of algebras is a monoidal equivalence.
For this reason, neither situation in Appendices~\ref{S:comonoid} and~\ref{S:counital}
is more general than the other.
\end{rem}

\medskip

The following is the analog of Proposition~\ref{P:cirmod} for augmented algebras.

\bp\label{P:augmod} Let $A$ and $B$ be augmented algebras. Let $M$ be a right
\odash{A}{module} and
$N$ a left \odash{B}{module}. Then $A\ten N$ is a left \odash{A\teninf
B}{module} via
\[(a\ten b)\cdot(a'\ten n):=\eta_B(b)aa'\ten n+\eta_A(a')a\ten bn\]
and $M\ten B$ is a right \odash{A\teninf
B}{module} via
\[(m\ten b')\cdot(a\ten b):=\eta_A(a)m\ten b'b+\eta_B(b')ma\ten b\,.\]
\ep
\bpf Similar to the proof of Proposition \ref{P:tensorprod}.
\epf

One may similarly show that the map
\[A\ten N\to A\cirprod N,\ \ a\ten n\mapsto \eta_A(a)n+a\ten n\]
is a morphism of left \odash{A\teninf B}{modules}, where $A\cirprod N$ is viewed
as a left \odash{A\teninf B}{module} by restriction via the morphism of
algebras~\eqref{E:infcir}.

Finally, we discuss the analog of the construction~\eqref{E:ciraction} for counital
$\infbis$, and apply these general considerations to the construction of an
$\infmod$. 

Let $A$ be a counital $\infbi$ and $N$ a left
\odash{A}{module}. It is possible to define a left \odash{A}{module} structure on
$A\ten N$, by restricting the structure of Proposition
\ref{P:augmod} along the morphism of augmented algebras $\Delta:A\to A\teninf A$.
By counitality, the action of $A$ on $A\ten N$ reduces to
\begin{equation}\label{E:action}
a\cdot(a'\ten n)=aa'\ten n+\eta(a')a_1\ten a_2n\,.
\end{equation}
We denote this module structure on the space $A\ten N$ by $A\teninf N$.

Our next result describes $\infmods$ over counital $\infbis$ in a way that
is analogous to the definition of Hopf modules over an ordinary Hopf algebra. 

Recall that a left Hopf module over a Hopf algebra $H$ is
a space $M$ that is both a left module and comodule over $H$ and for which the
comodule structure map $M\to H\ten M$ is a morphism of left
\odash{H}{modules}~\cite[Definition 1.9.1]{mon}. It is understood that $H\ten M$ is a left
\odash{H}{module} by restriction via the comultiplication of $H$. 

\bp\label{P:infmod} Let $A$ be a counital $\infbi$. Let $\lambda:A\ten N\to N$ be a left
\odash{A}{module} structure on $N$ and $\Lambda:N\to A\ten N$ a counital comodule structure
on $N$. Then $(N,\lambda,\Lambda)$ is an $\infmod$ over $A$ if and only if
$\Lambda:N\to A\teninf N$ is a morphism of left \odash{A}{modules}.
\ep
\bpf Write $\lambda(a\ten n)=an$ and $\Lambda(n)=n_{-1}\ten n_0$. According to 
\eqref{E:action}, 
\begin{align*}
a\cdot\Lambda(n) &=an_{-1}\ten n_0+\eta_A(n_{-1})\,a_1\ten a_2n_0\\
&=an_{-1}\ten n_0+a_1\ten a_2n\,,
\end{align*}
by counitality for $N$. Thus, 
$\Lambda$ is a morphism of \odash{A}{modules} if and only if
\[\Lambda(an)=an_{-1}\ten n_0+a_1\ten a_2n\,,\]
which is Axiom~\eqref{E:infmod} in the definition of $\infmod$.
\epf

Next, we make use of Proposition~\ref{P:infmod} to obtain 
the general construction of $\infmods$ of Example~\ref{Ex:nontrivial}.\!~3.

First, note that the tensor product construction of Proposition~\ref{P:augmod} is
associative, in the sense that if $A$, $B$ and $C$ are augmented algebras and
$N$ is a left \odash{C}{module}, then
\[(A\teninf B)\teninf N\cong A\teninf (B\teninf N)\]
as left \odash{A\teninf B\teninf C}{modules}. In fact, one has
\begin{multline*}
\bigl((a\ten b)\ten c\bigr)\cdot\bigl((a'\ten b')\ten n\bigr)=\\
\eta_C(c)\eta_B(b)aa'\ten b'\ten n+
\eta_C(c)\eta_A(a')a\ten bb'\ten n+
\eta_A(a')\eta_B(b')a\ten b\ten cn\\
=\bigl(a\ten (b\ten c)\bigr)\cdot\bigl(a'\ten (b'\ten n)\bigr)\,.
\end{multline*}

On the other hand, if $f:A\to B$ is a morphism of augmented algebras and $N$ is
a left \odash{C}{module}, then $f\ten\id_N:A\teninf N\to B\teninf N$ is a morphism of
left \odash{A\teninf C}{modules}, where $B\teninf N$ is an \odash{A\teninf C}{module} 
by restriction via the morphism of algebras $f\ten\id_C:A\teninf C\to B\teninf C$.

Now let us apply these considerations to a left module $N$ over a counital $\infbi$
$A$, $B=A\teninf A$, $C=A$ and $f=\Delta$. We obtain that 
\[\Delta\ten\id_N:A\teninf N\to A\teninf A\teninf N\]
is a morphism of left \odash{A\teninf A}{modules}. Hence, it is also a morphism of left
\odash{A}{modules}, by restriction via $\Delta$. An application of
Proposition~\ref{P:infmod} then yields the following

\bc\label{C:nontrivial} Let $A$ be a counital $\infbi$ and $N$ a left \odash{N}{module}.
Let $M=A\teninf N$, an \odash{A}{module} as in \eqref{E:action}. Define $\Lambda:M\to
A\ten M$ by 
\[\Lambda(a\ten n)=a_1\ten a_2\ten n\,.\]
With these module and comodule structures, $M$ is a left $\infmod$ over $A$.
\ec
\medskip
In this paper, quasitriangular $\infbis$ play an important role (Section~\ref{S:dendri}).
Our last result shows that the classes of counital $\infbis$ and
quasitriangular $\infbis$ are disjoint.

\bp\label{P:counital-quasi} If a quasitriangular $\infbi$ A is counital then $A=0$.
\ep
\bpf Let $r=\sum u_i\ten v_i$ be the canonical element and $\eta$ the counit.
According to~\eqref{E:principal}, we have 
\[\Delta(a)=\sum_i  u_i\ten v_ia-au_i\ten v_i\,.\]
 Applying $\eta\ten\id$
and using Lemma~\ref{L:counital} we deduce
\[a=\sum_i \eta(u_i)v_ia \text{ for every }a\in A\,.\]
Similarly, applying  $\id\ten\eta$ we deduce
\[a=-\sum_i au_i\eta(v_i) \text{ for every }a\in A\,.\]
Thus, $A$ has a left unit and a right unit. These must therefore coincide and
$A$ must be unital. But an $\infbi$ $A$ that is both unital and counital must be $0$
~\cite[Remark 2.2]{cor}.
\epf

\providecommand{\MR}{\relax\ifhmode\unskip\space\fi MR }
% \MRhref is called by the amsart/book/proc definition of \MR.
\providecommand{\MRhref}[2]{%
  \href{http://www.ams.org/mathscinet-getitem?mr=#1}{#2} }

\end{document}